\documentclass[12pt]{article}
\usepackage[utf8]{inputenc}

\usepackage{amsmath}
\usepackage{amscd}
\usepackage{amssymb}
\usepackage{enumerate}
\usepackage{indentfirst}
\usepackage{latexsym}
\usepackage{multicol}
\usepackage{verbatim}

\newcommand{\Int}{\displaystyle \int}
\newcommand{\Frac}{\displaystyle \frac}

\usepackage{theorem}
\numberwithin{equation}{section}
\def\Id{\hbox{\rm Id}}
\def\div{\hbox{\rm div}\,  }

\def\divA{\, \hbox{\rm div}_A\,  }
\def\divu{\, \hbox{\rm div}_u\,  }
\def\divx{\, \hbox{\rm div}_x\,  }
\def\divy{\, \hbox{\rm div}_y\,  }

\def\Deltau_{\Delta_u}

\newcommand{\R}{{\mathbb R}}
\newcommand{\N}{{\mathbb N}}

\def\d{\partial}

\def\ep{\varepsilon}

\def\cB{{\mathcal B}}

\let\tilde=\widetilde

\newcommand{\du}{\delta\! u}

\newcommand{\dA}{\delta\! A}

\newcommand\dP{\delta\!P}

\newtheorem{lem}{Lemma}

\newtheorem{prop}{Proposition}
\newtheorem{theo}{Theorem}
\newtheorem{rem}{Remark}

\newenvironment{p}{
\smallbreak\noindent\textit{\textbf{Proof:}}~}
{\hfill\rule{2.1mm}{2.1mm}\smallbreak
}

\begin{document}

\centerline{\large \bf Incompressible flows with piecewise constant density}

\begin{center}

{Rapha\"el Danchin$^1$ and Piotr Bogus\l aw Mucha$^2$}

\medskip 

{1. Universit\'e Paris-Est, LAMA, UMR 8050 and Institut Universitaire de France,}

{ 61 avenue du G\'en\'eral de Gaulle,
94010 Cr\'eteil Cedex, France.}

{E-mail: danchin@univ-paris12.fr}

{2. Instytut Matematyki Stosowanej i Mechaniki,
 Uniwersytet Warszawski, }

{ul. Banacha 2,  02-097 Warszawa, Poland.} 

{E-mail: p.mucha@mimuw.edu.pl}

\end{center}

\date\today
\begin{abstract} We investigate the incompressible Navier-Stokes equations  with variable density. 
The aim is to  prove existence and uniqueness results in the case of  \emph{discontinuous} initial density.
In dimension $n=2,3,$  assuming only that the initial density is bounded and bounded away from zero, and that the initial velocity is smooth enough, we get the local-in-time existence of unique solutions. Uniqueness holds in any dimension and for a wider class of 
 velocity fields. Let us emphasize that  all those  results are true for  piecewise constant densities with arbitrarily large jumps.
Global results are established in dimension two if the density is close enough to a positive constant, 
and in $n$-dimension  if, in addition,  the initial velocity is small. 
The Lagrangian formulation for describing the flow plays a key role in the analysis that is proposed in the present paper. 
\end{abstract}

\noindent
{\it MSC:} 35Q30, 76D05

\noindent
{\it Key words:} Inhomogeneous Navier-Stokes equations, critical regularity,  piecewise constant density, large jumps, Besov spaces, Lagrangian coordinates,  discontinuous data.

\section*{Introduction}

Incompressible flows are often modeled by the \emph{homogeneous} Navier-Stokes equations~: 
that is  the density of the fluid  is assumed to be a constant. However in many applications as blood flows or models of rivers, 
although the fluid is practically incompressible, the density can not be considered
as a constant quantity, as a consequence of the complex structure of the flow due to  e.g.  a mixture of fluids or
 pollution (see e.g.  \cite{AS,Ar,Ku,LZWL,San}). 
This makes us look at the density as  a nonnegative unknown function which has constant values along the stream lines. The simplest model which can capture such a physical property
is the  so-called \emph{inhomogeneous Navier-Stokes system}:
\begin{equation}\label{NSE}
\begin{array}{lcr}
\rho_t +v \cdot \nabla \rho =0 & \mbox{in} & \Omega \times (0,T), \\
\rho v_t + \rho v \cdot \nabla v - \nu  \Delta v + \nabla Q =0& \mbox{in} & \Omega \times (0,T), \\
\div v =0& \mbox{in} & \Omega \times (0,T), \\
v=0 & \mbox{on} & \d\Omega \times (0,T), \\
v|_{t=0}=v_0, \qquad \rho|_{t=0}=\rho_0 & \mbox{in} & \Omega.
\end{array}
\end{equation}
The unknown functions are: $\rho$ -- the density of the fluid, $v$ -- its velocity field and $Q$ -- its pressure. The constant positive viscosity coefficient is denoted by $\nu$.
We consider the  cases where  $\Omega$ is a bounded domain of $\R^n,$ or the whole space $\R^n,$
and we focus mainly on the physically relevant space dimensions $n=2,3.$ 

%

The goal of the present paper is to revisit  results concerning the well-posedness issue of \eqref{NSE}. We 
concentrate our analysis on the regularity of density. In our recent work \cite{DM},  we established 
the existence and uniqueness of solutions to \eqref{NSE} in a critical regularity framework 
which allowed the initial density to be discontinuous.  However,  a smallness condition over 
the jumps was required there. 
In the present work, we want to discard this smallness condition.
At the same time, to simplify the presentation, we do not
strive for optimal assumptions as concerns the velocity and  assume
the viscosity coefficient $\nu$ to be constant. 

Let us recall (see in particular \cite{Desjardins} and  the textbooks \cite{AKM,Lions}) that, roughly, from the qualitative viewpoint the classical 
results for the homogeneous Navier-Stokes equations carry out to \eqref{NSE} :
on the one hand global (possibly non unique) weak solutions with finite energy may be constructed 
and on the other hand,  if the density is smooth enough, bounded and bounded away from zero,  
then   global-in-time existence and uniqueness results are available in dimension two 
for arbitrarily large data, and if the velocity is small in dimension three. 
These latter results require  relatively high regularity of the density, though.  In particular it has 
to be at least continuous, and to have some fractional derivatives in suitable Lebesgue spaces 
(see e.g.  \cite{LS} or \cite{D1}). It is worthwhile to emphasize that for smooth densities one may show the existence of unique solutions
even for vacuum states \cite{CK} provided the initial data satisfy some compatibility condition.
{}From the viewpoint of applications such results are not so  satisfactory: we wish to consider fluids with 
e.g.   piecewise constant densities, a pattern which is of interest to  model a mixture of two fluids.
\medbreak

The results of the paper are split into two groups:
\begin{itemize}
\item
The first group concerns  uniqueness and  local-in-time existence results in the case 
where  the initial density is just an $L_\infty$ positive function bounded away from zero. In particular, one may  consider   
piecewise constant densities \emph{with arbitrary large  jumps.} As regards  the  existence issue,  we have to restrict ourselves
to the (physically relevant) dimensions  $n=2,3$. 
\item The second group concerns the global-in-time existence issue. 
Here we have  to make a smallness assumption over 
the density which, in the case of piecewise constant initial density, implies that  the jumps have to be small.
Assuming enough smoothness over the velocity, this enables us to prove global existence for (possibly) large
velocity if $n=2,$ and for small velocity if $n \geq 3$ (an assumption which is also  required for the homogeneous Navier-Stokes equations, anyway).
\end{itemize}

As explained above, in the present paper, we aim at doing minimal assumptions over the density but
 we do not strive for optimal regularity of the velocity function. As our method relies  on  estimates for 
the Stokes system with merely bounded coefficients, the (rather high) regularity of the velocity is somehow prescribed by the technique. 
An approach  to the issue 
of sharp regularity has been done in \cite{DM} in the critical Besov spaces
setting. However in \cite{DM} we were able to capture discontinuous  density with small jumps only.

The rest of the paper unfolds as follows. The main results are presented in the first section.
Then, some preliminary estimates involving the evolutionary Stokes system are proved. 
Section \ref{s:lagrangian} is devoted to the derivation of System \eqref{NSE} in Lagrangian coordinates. 
In Section \ref{s:uniqueness}, we concentrate on the proof of uniqueness results
whereas existence results are proved in the last two sections.  Technical estimates involving the divergence equation are presented in the Appendix.


\section{Main results}\label{s:main}

Let us first recall the basic energy equality for  System \eqref{NSE}
which may be (formally) derived  by testing $(\ref{NSE})_2$ by $v$: 

\begin{lem}\label{l:ene}
 Let $(\rho,v)$ be a sufficiently smooth solution to \eqref{NSE} over $\Omega\times[0,T].$ Then there holds
\begin{equation}\label{energy}
\int_\Omega (\rho|v|^2)(t,x)\,dx+2\nu\! \int_0^t\!\!\int_\Omega|\nabla v(\tau,x)|^2\,dx\,d\tau=\int_\Omega (\rho|v|^2)(0,x)\,dx\!\!\quad\hbox{for all }
\ t\in[0,T].
\end{equation}
\end{lem}
Subsequently if $\rho_0$ is positive and bounded away from zero  
and $v_0$ is in $L_2(\Omega)$ then we  get a control over
$v$ in $L_\infty(0,T;L_2(\Omega))$  and $\nabla v$ in $L_2(\Omega \times (0,T)).$
Under very rough regularity assumptions (much less than assumed here), 
the (formal) energy equality \eqref{energy}  provides us with an  information about low norms of the velocity, which turns out to be
crucial for the proof of global results (see in particular the monograph by \cite{Lions} and the references therein,
as regards the proof of  global weak solutions with finite energy). 
Note  that \eqref{energy}  gives some regularity information 
over the velocity even for very rough density. We shall see further in the paper
a way to get even more regularity information over the velocity without assuming more on the density.

\smallskip
Before listing the main results of the paper, let us introduce a few notation. 
Concerning the derivatives of functions $f$ depending on both
the time variable $t$ and the space variable $x,$ we denote by $f_t$ the time derivative and
by $Df$ the Jacobian matrix of $f$ with respect to the space variable, namely
$(Df)_{i,j}=\d_j f^i.$ 
The notation $\nabla f$ is reserved for ${}^T\!(Df).$

The Lebesgue spaces  
of measurable functions  with integrable  $p$-th power is denoted by $L_p(\Omega)$.
More generally, if $m\in\N$ then $W_p^m(\Omega)$ denotes the 
set of $L_p(\Omega)$ functions with derivatives of order less than or equal to $m$ in $L_p(\Omega).$
 Since the Navier-Stokes equations are of  parabolic type, it is also natural 
to introduce  \emph{parabolic} Sobolev spaces  $W^{2,1}_{q,p}(\Omega \times (0,T))$ that is  the closure 
 of smooth functions for the norm
\begin{equation}
\|u\|_{W^{2,1}_{q,p}(\Omega \times (0,T))} = \|u,\d_t u\|_{L_p(0,T;L_q(\Omega))}
+\|u\|_{L_p(0,T;W^2_p(\Omega))}.
\end{equation}
Granted with parabolic spaces, one may now define  Besov spaces over $\Omega$ as
the following trace space:
\begin{equation}\label{xx1}
 B^{2-2/p}_{q,p}(\Omega)= \Bigl\{ f:\Omega\rightarrow\R\ \hbox{measurable s.t.} \ f=\bar f|_{t=0}\!\ \hbox{ for some }\  \!\bar f \!\in\! W^{2,1}_{q,p}(\Omega \times (0,1))\Bigr\}\cdotp
\end{equation}
The norm can be defined from the above definition  as a suitable infimum (for  more details concerning the Besov spaces we refer to \cite{BIN,Triebel}).
\smallbreak

Our first result  states the uniqueness of  solutions with merely bounded density, provided the initial velocity 
is  smooth enough. 
\begin{theo}\label{th:uniq} Let $n \geq 2$. 
Assume  that $\Omega$ is $\R^n$ or a $C^2$ bounded domain of $\R^n.$
 Let $(\rho^1,v^1,Q^1)$ and $(\rho^2,v^2,Q^2)$ be two solutions to \eqref{NSE} 
with the same initial data, and density bounded and bounded away from $0.$ 
Suppose moreover that for $k=1,2,$
\begin{itemize}
\item Case $n=2$:  there exists $q>2$ such that  $v^k \in W^{2,1}_{q,2}(\Omega\! \times\! (0,T))$ and  $\nabla Q^k \in L_2(0,T;L_q(\Omega))$,
\item Case $n \geq 3$:     $v^k \in W^{2,1}_{n,2}(\Omega \times (0,T))$, $\nabla Q^k \in L_2(0,T;L_n(\Omega))$
and, in addition, $\nabla v^k\in L_2(0,T;L_\infty(\Omega)).$
\end{itemize}
Then $v^1\equiv v^2,$ $\nabla Q^1\equiv\nabla Q^2$ and $\rho^1\equiv \rho^2$. 
\end{theo}

\begin{rem} 

As regards the inhomogeneous incompressible Navier-Stokes equations, to our knowledge, the ``best'' uniqueness 
result with no smallness condition over the density is due to P. Germain in \cite{Germain}. 
It does not apply to  solutions with piecewise constant densities, though.
\end{rem}

The second result complements Theorem \ref{th:uniq}. It delivers existence of local-in-time regular and unique solutions
in dimensions $2$ and $3.$  Again, the initial density just has to be bounded and bounded away from vacuum.

\begin{theo}\label{th:exist} 
 Let $n=2,3$ and $\Omega$ be a $C^2$ bounded domain  or be $\R^n$. 
 Let   $\rho_0$ satisfy 
 \begin{equation}\label{eq:rho}
 m<\rho_0<M
 \end{equation}
 for some positive constants,  and  $v_0\in W^2_2(\Omega)$ be such that $\div v_0=0$ and 
$v_0|_{\d\Omega}=0$. Let  $n^*=2\bigl(\frac{n+2}n\bigr).$ 
There exists a unique solution $(\rho,v)$ to System \eqref{NSE} on a time interval $[0,T]$ for some $T>0$ such that
$\rho(t,\cdot)$ satisfies \eqref{eq:rho} for all $t\in[0,T]$ and 
$$
v \in W^{2,1}_{n^*,n^*}(\Omega \times (0,T)), \quad v_t\in L_\infty(0,T;L_2(\Omega))\mbox{ \ \ and \ \ } \nabla v_t \in L_2(\Omega \times (0,T)).
$$ 
\end{theo}

\begin{rem} The critical Sobolev embedding ensures that $W^2_2(\Omega)$ is continuously embedded
in  the Besov space $B^{2-2/n^*}_{n^*,n^*}(\Omega).$ Keeping in mind the definition of this space given in \eqref{xx1}, 
the appearance of the parabolic Sobolev space  $W^{2,1}_{n^*,n^*}(\Omega \times (0,T))$ in the above statement does not
come up as a surprise. The $W^2_2(\Omega)$ assumption for $v_0$ is needed to ensure that 
$(\d_tv+v\cdot\nabla v)|_{t=0}$ is in $L_2(\Omega).$ 
At the same time, owing to the low regularity of the density, 
we do not  know how to propagate the  $W^2_2(\Omega)$ regularity for the velocity.
\end{rem}

Proofs of Theorems \ref{th:uniq} and \ref{th:exist} are based on the analysis of \eqref{NSE} in the Lagrangian coordinates defined by the stream lines. Since the density is merely bounded
there is an obstacle to apply any bootstrap method in order to improve the regularity of the velocity. The main difficulty is located in the term with the time derivative. To obtain a better 
information about $v_t,$ we adopt  techniques from the compressible Navier-Stokes system \cite{Mu-CM,MZ2} (concerning uniqueness criteria for the compressible Navier-Stokes system in Lagrangian formulation, see also the recent work by D. Hoff).
Roughly speaking, we differentiate the (Lagrangian) velocity equation once with respect to time, then 
apply an energy method. 
This  approach via the Lagrangian coordinates requires only  $L_\infty$ bounds (by above and by below)  for  the density,  provided  the velocity has high regularity.   That the  density is time-independent in the Lagrangian setting, hence is just a given function, is of course fundamental.
In comparison, in \cite{D1} where the Eulerian framework is used,  the initial density has to be in the Besov space $B^{n/2}_{2,1}(\R^n)$ (which, roughly, means that it has $n/2$ derivatives in $L_2(\R^n)$)
but the initial   velocity therein has only \emph{critical regularity}, namely it is in $B^{n/2-1}_{2,1}(\R^n)$ (to be compared with 
$B^2_{2,1}(\R^n)$ and $n=2,3$ here).

 To highlight   consequences of Theorems \ref{th:uniq} and \ref{th:exist}, let us consider  the case where the initial
 divergence-free velocity field is in $W^2_2(\Omega)$ (and vanishes at the boundary), and the initial
  density $\rho_0$ is\footnote{Here 
$\chi_{A_0}$ stands for  the characteristic function of the set $A_0$.} 
\begin{equation}\label{rem1}
\rho_0=m+\sigma \chi_{A_0},
\end{equation}
where $m,\sigma$ are positive constants and $A_0$ is a set with a $C^1$ boundary. 
The velocity field $v$  given  by Theorem \ref{th:exist} is Lipschitz with respect to the space
variable  hence generates a unique $C^1$ flow $X$ defined by
$$
X(t,y)=y+\int_0^tv(\tau,X(\tau,y))\,d\tau.
$$
Therefore, the density at time $t$ is given by 
\begin{equation}\label{rem2}
\rho(t,\cdot)=m+\sigma \chi_{A(t)},
\quad\hbox{with }\  A(t):=X(t,A_0).
\end{equation}
As  the flow $X$  is at least $C^1$,   the initial regularity of the boundary 
of $A(t)$ is preserved and   any geometrical catastrophe (e.g. breaking down  or self-intersections 
of the boundary)  will not appear : if $A_0$ is diffeomorphic to a ball, then $A(t)$ is diffeomorphic to a ball, too.
The above case shows that the system \eqref{NSE} can model an interaction of two fluids separated by a free interface. Although tracking the regularity of the boundary $\d A(t)$ is not the main topic of this paper, we see that  Theorem \ref{th:exist} ensures
that the  $C^1$  or $C^{1,\alpha}$ regularity (with  $\alpha$ small enough) of $\d A(t)$ is preserved during the evolution. 
In other words, we have  partially solved in an indirect way a complex free boundary
problem which has been left as an open question by P.-L. Lions in his book \cite{Lions}. 
Let us emphasize that the  standard  approach for solving problems of such type  requires very technical considerations 
(see e.g. \cite{Ables,Shimizu,TW}).
Furthermore,   with our approach, there are no requirements  concerning the regularity of the boundary of the set $A_0$:
our results hold for any measurable set $A_0$. 

\medbreak

The above results concern local-in-time analysis. In order to obtain global-in-time  solutions,  we have  to 
assume that the jumps of the initial density are small enough. 
The following theorem states that under this sole assumption over the density, and for sufficiently smooth (possibly large) initial velocity fields, 
global existence holds true.
\begin{theo}\label{th:lar}
Let $\Omega$ be a $C^2$ bounded two-dimensional set, or  be $\R^2$. 
There exists a constant $c$ depending only on $\Omega$ and such that if
$\rho_0 \in L_\infty(\Omega)$ satisfies
\begin{equation}\label{den-str}
 \frac{\displaystyle \sup_{x \in \Omega} \rho_0(x) - \inf_{x\in \Omega} \rho_0(x)}{\displaystyle \inf_{x\in \Omega} \rho_0(x)} \leq c
\end{equation}
then for all  $v_0 \in B^{1}_{4,2}(\Omega)\cap L_2(\Omega)$ with $\div v_0=0$ and $v_0|_{\d\Omega}=0,$
 there exists a unique global-in-time solution to System \eqref{NSE} such that \eqref{energy} is satisfied and that, for all $T>0,$
\begin{equation}\label{l2}
v \in W^{2,1}_{4,2}(\Omega \times (0,T)),\ \  \nabla Q\in L_2(0,T;L_4(\Omega)) \mbox{ \ and \  } \rho \in L_\infty(\Omega \times (0,T)).
\end{equation}
\end{theo}

In dimension $n\geq3,$ getting global-in-time strong solutions requires also the initial velocity to be small, 
(an assumption which is  needed in the homogeneous case, anyway). Here is our statement:
\begin{theo}\label{th:bdd}
Let $\Omega$ be a bounded $n$-dimensional $C^2$ domain. Let $\rho_0 \in L_\infty(\Omega)$
be positive and bounded away from $0,$ and
$v_0 \in B^{2-\frac 2q}_{q,p}(\Omega)$ with $1<p<\infty,$  $n<q<\infty$ and  $2-2/p\not=1/q.$ 
There exist two constants $c$ and $c'$ depending only on $\Omega,$ $p$ and $q$ and such that if
\begin{equation}\label{den-str-1}
 \frac{\displaystyle \sup_{x \in \Omega} \rho_0(x) - \inf_{x\in \Omega} \rho_0(x)}{\displaystyle \inf_{x\in \Omega} \rho_0(x)} < c
 \quad\hbox{and}\quad
 \|v_0\|_{B^{2-\frac 2p}_{q,p}(\Omega)} \leq c'\nu,
\end{equation}
then there exists  a unique global-in-time solution  to the inhomogeneous Navier-Stokes system \eqref{NSE} such that 
$$
v \in W^{2,1}_{q,p}(\Omega \times\R_+),\ \  \nabla Q\in L_p(\R_+;L_q(\Omega)) \mbox{ \ and \  } \rho \in L_\infty(\Omega \times \R_+).
$$
Furthermore, there exist two positive constants $\alpha$ and $C$ depending only on $\Omega,$ $p,$ $q$ and of
the lower and upper bounds for $\rho_0$ so that for all $t>0,$
$$
\|v\|_{W^{2,1}_{q,p}(\Omega \times(t,t+1))}+\|\nabla P\|_{L_p(t,t+1;L_q(\Omega))}\leq Ce^{-\alpha t} \|v_0\|_{B^{2-\frac 2p}_{q,p}(\Omega)}.
$$
\end{theo}

 Theorems \ref{th:lar} and \ref{th:bdd} follow from classical maximal regularity techniques.
  The smallness conditions \eqref{den-str} and \eqref{den-str-1} allow 
to treat the  oscillations of the density  as a perturbation that may be put  in  the right-hand side  of the estimates.

At the end we would like to underline that most of our results hold for bounded domains and $\R^n$. 
The case of the whole space is easier: there is no boundary condition and  solving  the divergence equation is simpler, too. 
One exception is Theorem \ref{th:bdd}  where the 
boundedness of the domain is essential here as it provides exponential decay of the energy norm
 (the whole space case is tractable  under stronger conditions over the density, 
see our recent work in \cite{DM}).


\section{Some linear estimates}\label{s:linear}

A fundamental role will be played by the Stokes system, that is the following 
linearization of the velocity  equation in \eqref{NSE}:
\begin{equation}\label{stokes}
\begin{array}{lcr}
mu_t - \nu \Delta u + \nabla Q= f & \mbox{in} & \Omega \times (0,T), \\[1ex]
\div u= \div R & \mbox{in} & \Omega \times (0,T), \\[1ex]
u=0 & \mbox{in} &\d \Omega \times (0,T), \\[1ex]
u|_{t=0}=u_0 & \mbox{in} & \Omega,
\end{array}
\end{equation}
where $m$ and $\nu$ are positive constants.
\smallbreak
We shall make an extensive use of 
the following solvability result for the Stokes system
in the $L_p(0,T; L_q(\Omega))$ framework:

\begin{theo}\label{th:stokes} Let $\Omega$ be a $C^2$ bounded subset of  $\R^n.$ Let 
 $1<p,q<\infty,$ $u_0\in B^{2-\frac 2p}_{q,p}(\Omega),$ $f\in L_p(0,T;L_q(\Omega))$, $R\in W^1_p(0,T;L_q(\Omega))$
so that  
$\div R \in L_p(0,T;W^1_q(\Omega)).$  
Suppose that
$$\div u_0=\div R|_{t=0}\quad\hbox{and}\quad\vec n \cdot R|_{\d \Omega \times (0,T)} =0.$$
If $2-2/p> 1/q,$  assume in addition that $u_0=0$ at the boundary, otherwise
we assume only $u_0\cdot \vec n =0$ at $\d\Omega$\footnote{For simplicity we exclude the case $2-2/p=1/q$.}. 
Then there exists a unique solution to \eqref{stokes} such that $u\in W^{2,1}_{q,p}(\Omega \times (0,T))$, $\nabla P \in L_p(0,T;L_q(\Omega)),$ and the following estimate is valid:
\begin{multline}\label{est-sto}
\|mu_t, \nu \nabla^2 u, \nabla P\|_{L_p(0,T;L_q(\Omega))}+ m^{\frac1p}\nu^{\frac1{p'}}\sup_{0\leq t \leq T}\|u(t)\|_{B^{2-\frac 2p}_{q,p}(\Omega)}  \\ 
\leq C 
(\|f,m R_t\|_{L_p(0,T;L_q(\Omega))} + \|\nu \div R\|_{L_p(0,T;W^1_q(\Omega))} + m^{\frac1p}\nu^{\frac1{p'}}\|u_0\|_{B^{2-\frac 2p}_{q,p}(\Omega)}),
\end{multline}
where $C$ is independent of $m$, $\nu$ and $T$. 
\end{theo}

\begin{p}
In the case $R\equiv0,$ this result is classical (see e.g. 
 \cite{GS,MS} and the appendix of \cite{D}). 
The general case follows from this particular case  once constructed a suitable vector-field $w:\Omega\times (0,T) \to
\R^n$ fulfilling
\begin{equation}\label{eqdiv}
 \div w=\div R \mbox{~~in }\Omega, \qquad w=0 \mbox{~~ at } \d \Omega.
\end{equation}
Taking for granted the existence of such a vector-field, the strategy is simple : 
we set $v=u-w$ and we gather that $v$ has to satisfy 
$$\begin{array}{lcr}
mv_t - \nu \Delta v + \nabla Q= f-mw_t+\nu\Delta w & \mbox{in} & \Omega \times (0,T), \\[1ex]
\div v= 0 & \mbox{in} & \Omega \times (0,T), \\[1ex]
v=0 & \mbox{in} &\d \Omega \times (0,T), \\[1ex]
v|_{t=0}=u_0-w_0 & \mbox{in} & \Omega.
\end{array}
$$
Therefore, in order to reduce our study to the case $R\equiv0,$ the vector-field $w$ is required to 
satisfy $w_t,D^2w\in L_p(0,T;L_q(\Omega))$ (note that this will imply that 
$w_0\in B^{2-2/p}_{q,p}(\Omega),$ see \eqref{xx1}).  
The fact that such a solution to \eqref{eqdiv} does exist 
is granted by the following lemma, the proof of which is postponed in Appendix (see Proposition \ref{p:bog}):
\end{p}

\begin{lem}\label{l:bog} 
 Let $R(t,\cdot)$ be a family of vector-fields defined over the $C^2$ bounded domain $\Omega,$ parameterized by $t\in (0,T).$
 Assume that, for some $1<q<\infty$ and $1\leq p\leq\infty$  we have
  $\div R \in L_p(0,T;W^1_q(\Omega))$, $R,R_t \in L_q(0,T;L_p(\Omega))$ and $R\cdot \vec n =0$ at the boundary. 

Then there exists a  vector-field  $w$ in $L_p(0,T;W^2_q(\Omega))$ vanishing on $\d\Omega,$  fulfilling 
$$
\div w=\div R\quad\hbox{and}\quad \div w_t=\div R_t\quad\hbox{in }\ \Omega
 $$
 and  the following estimates: 
\begin{eqnarray}\label{eq:bog1}
 &\|w\|_{L_p(0,T;W^2_q(\Omega))} \leq C\|\div R\|_{L_p(0,T;W^1_q(\Omega))},\\ \label{eq:bog2}
&\|w_t\|_{L_p(0,T;L_q(\Omega))} \leq C\|R_t\|_{L_p(0,T;L_q(\Omega))}
\end{eqnarray}
for some constant $C$ depending only on $q$ and $\Omega.$ 
\end{lem}

\begin{rem}\label{r:stokes} The whole space case is easier to deal with for we do not have to take care of boundary conditions (apart
from suitable decay at infinity given by the functional setting). Indeed, in  order to solve \eqref{eqdiv}, one may
set 
$$
w=-\nabla(-\Delta)^{-1}\div R. 
$$
As the corresponding Fourier multiplier is homogenous of degree $0,$ we  readily get \eqref{eq:bog1} and \eqref{eq:bog2}. 
Therefore, arguing as above and using the standard maximal regularity result for the Stokes system in $\R^n,$ 
we conclude to 
 Theorem \ref{th:stokes}  in the case $\Omega=\R^n$ if  the Besov space $B_{q,p}^{2-\frac2p}(\Omega)$
 is  replaced with the \emph{homogeneous Besov space} $\dot B^{2-\frac2p}_{q,p}(\R^n)$ and
 $W^k_q(\Omega),$ by its homogeneous version $\dot W^k_q(\R^n).$
 \end{rem}

Theorem \ref{th:stokes} can be viewed as a classical result. In order to prove Theorem \ref{th:exist} we need 
to adapt it to the case of \emph{variable} coefficients. Below, we focus 
on the $L_2$ case where only the boundedness of coefficients is needed. 
\begin{lem}\label{l:l2stokes} Let $\Omega$ be a bounded domain of $\R^n,$ or $\R^n.$
 Let $\eta \in L_\infty(\Omega)$ be a time independent positive function, bounded away from zero, and
 $R$ satisfy the above boundary conditions. Then the solution $(u,\nabla P)$ with $u\in W_{2,2}^{2,1}(\Omega\times(0,T))$
 and $\nabla P\in L_2(\Omega\times(0,T))$ to  
 the system 
\begin{equation}\label{u2}
 \begin{array}{lcr}
  \eta u_t -\nu \Delta u + \nabla P =f & \hbox{ in } & \Omega\times (0,T),\\
\div u = \div R & \mbox{ in } & \Omega\times (0,T),\\
u=0 & \mbox{ on } & \d \Omega\times (0,T),\\
u|_{t=0}=u_0 & \mbox{ on } & \Omega,
 \end{array}
\end{equation}
fulfills 
\begin{multline}\label{u3}
\sqrt\nu\sup_{0\leq t \leq T} \|\nabla u(t)\|_{L_2(\Omega)}+
  \|u_t,\nu\nabla^2u,\nabla P\|_{L_2(\Omega \times(0,T))}
 \leq C\big( \|f,R_t\|_{L_2(\Omega \times (0,T))} \\
+\nu\|\div R\|_{L_2(0,T;W^1_2(\Omega))} +\sqrt\nu\|\nabla u_0\|_{L_2(\Omega)}\big), 
\end{multline}
where  $C$ depends on $\inf\eta$ and $\sup\eta,$  but is independent of $T$ and $\nu.$
\end{lem}

\begin{p}
 First we remove the right-hand side  of $\eqref{u2}_2$ by means of Lemma \ref{l:bog} (or the remark that follows if $\Omega=\R^n$): 
 we introduce a vector-field  $w$ fulfilling \eqref{eqdiv} such that $w\in W^{2,1}_{2,2}(\Omega \times(0,T))$ 
with the following bound
\begin{equation}\label{u5}
 \|\d_tw,\nu D^2w\|_{L_2(0,T;L_2(\Omega))} \leq C(\|\div R\|_{L_2(0,T;W^1_2(\Omega))} + \|R_t\|_{L_2(\Omega \times (0,T))}).
\end{equation}
Hence we may reduce the proof  to  the case  $R \equiv 0$. 
Now,  we observe that testing by $u_t$ gives:
$$
 \int_\Omega \eta |u_t|^2\, dx +  \frac\nu2\frac{d}{dt}  \int_\Omega |\nabla u|^2\, dx
 =\int_\Omega f\cdot u_t\,dx.
$$
Therefore, integrating in time yields
\begin{equation}\label{u6}
\nu\|\nabla u(t)\|_{L_2(\Omega)}^2+\int_0^t\|\sqrt\eta\,u_t\|_{L_2(\Omega)}^2\,d\tau\leq
\nu\|\nabla u_0\|_{L_2(\Omega)}^2+\int_0^t\|f/\sqrt\eta\|_{L_2(\Omega)}^2\,d\tau.
\end{equation}
Since $\eta$ is a positive time independent function which is pointwise bounded from below and above, we obtain 
\begin{equation}\label{u7}
 \|u_t\|_{L_2(\Omega \times (0,T))}+\sup_{0\leq t \leq T} \sqrt\nu\|Du(t)\|_{L_2(\Omega)} \leq C(\|f\|_{L_2(\Omega \times (0,T))} 
  +\sqrt\nu\|Du_0\|_{L_2(\Omega)}).
\end{equation}
In order to estimate  $D^2u$ and $DP,$  we rewrite  \eqref{u2} as
\begin{equation}\label{u8}
 \begin{array}{lcr}
   -\nu \Delta u + \nabla P =f - \eta u_t \qquad \qquad & \mbox{ in } & \Omega\times (0,T),\\
\div u = 0 & \mbox{ in } & \Omega\times (0,T),\\
u=0 & \mbox{ at } & \d \Omega\times (0,T).
 \end{array}
\end{equation}

If $\Omega$ is a $C^2$ bounded domain then the solvability of \eqref{u8} in the $L_2$ framework is clear (see e.g. \cite{Galdi}, Th. 6.1, page 231), thus taking into account bounds \eqref{u5}  and \eqref{u7} we get \eqref{u3}.  Lemma \ref{l:l2stokes} is  proved.
In the $\R^n$ case, one may just notice that
$\nabla P=-\nabla(-\Delta)^{-1}\div (f-\eta u_t).$ As $f-u_t$ is in $L_2(\Omega\times(0,T)),$ we still 
get the result, first for $\nabla P,$ and next for $\nabla^2u.$
\end{p}


\section{The Lagrangian coordinates}\label{s:lagrangian}

A fundamental point of our analysis is the use of  Lagrangian coordinates. In order to define them we solve the following 
ordinary differential equation (treating $y$ as a parameter):
\begin{equation}\label{i2}
\frac{d X(t,y)}{dt} = v( t,X(t,y)), \qquad X(t,y)|_{t=0}=y.
\end{equation}
This  leads to the following  relation between the Eulerian coordinates $x$ and the Lagrangian coordinates $y$:
\begin{equation}\label{i3}
X(t,y)=y+\int_{0}^t v(\tau,X(\tau,y))\, d\tau.
\end{equation}

Let us list a few basic properties for the Lagrangian change of variables:

\begin{prop}\label{p:lag}
 Suppose that  $v\in L_1(0,T;W^1_\infty(\Omega))$ with $v \cdot \vec n|_{\d \Omega}=0.$ Then the solution to System \eqref{i2} exists on the time interval $[0,T]$, $X(t,\Omega)=\Omega$ for all $t\in[0,T),$ and  $D_yX\in L_\infty(0,T;L_\infty(\Omega))$ with in addition
\begin{equation}\label{i5}
 \|D_yX(t)\|_{L_\infty(\Omega)}\leq \exp\biggl(\int_0^t\|D_xv\|_{L_\infty(\Omega)}\,d\tau\biggr).
\end{equation}
Furthermore 
\begin{equation}\label{i5a}
 X(t,y)=y +\int^t_0 u(t',y)\,dt' \mbox{ \ \ with \ } u(t,y):=v(t,X(t,y))
\end{equation}
so that $DX$  satisfies
\begin{equation}\label{i6}
 D_yX(t,y)=\Id +\int^t_0 D_yu(t',y)\,dt'.
\end{equation}
Let $Y(t,\cdot)$ be the inverse diffeomorphism of $X(t,\cdot).$
Then 
\begin{equation}\label{i6b}
D_xY(t,x)=(D_yX(t,y))^{-1}\quad\hbox{with }\ x=X(t,y)
\end{equation}
and, if $\displaystyle \int_0^t |D_yu(t,y)|\,dt'\leq 1/2$ then 
\begin{equation}\label{i7}
 \left| D_xY(t,x) - \Id\right| \leq 2\int_0^t |D_yu(t',y)|\,dt'.
\end{equation}
Finally,  if $v\in L_1(0,T;W^2_\infty(\Omega))$ then $D_yX\in  L_\infty(0,T;W^1_\infty(\Omega))$ and 
\begin{equation}\label{i6a}
\left|D^2_yX(t,y)\right|\leq e^{\int_0^t|D_xv(t',X(t',y))|\,dt'}
\int_0^t|D^2_xv(t',X(t',y))|e^{\int_0^{t'}|D_xv(t'',X(t'',y))|\,dt''}dt'.
\end{equation}
and  if $v \in L_1(0,T;W^s_p(\Omega))$ with $s>\frac np +1$, then $D_yX-\Id\in L_\infty(0,T;W^{s-1}_p(\Omega))$.
\end{prop}
\begin{p}
 The existence of $X$ for $(t,y) \in (0,T)\times\Omega$ follows from
 the standard  ODE theory, a consequence of Picard's theorem.  
 Inequalities \eqref{i5} and \eqref{i6a} follow from \eqref{i3} by differentiation and Gronwall lemma. 
 The higher regularity stems from the fact that, under our
 assumptions,  $W^{s-1}_p$ is an algebra (the reader may refer to the appendix of \cite{DM} for the proof
 of similar results in a slightly different context).

Equation \eqref{i6} follows from \eqref{i5a}, by differentiation. Then  \eqref{i7} 
comes from \eqref{i6b} provided $D_yX - \Id$ is small enough:  indeed, we have 
$$
D_xY=(\Id+(D_yX-\Id))^{-1}=\sum_{k=0}^{+\infty}(-1)^k\biggl(\int_0^tD_yu(t',y)\,dt'\biggr)^k.
$$
This  yields \eqref{i6}.
\end{p}

\medbreak

Let us now derive the Navier-Stokes equations  \eqref{NSE}  in \emph{ the  Lagrangian coordinates}: we set 
\begin{equation}\label{i4}
\eta(t,y):=\rho(t,X(t,y)), \quad u(t,y):=v(t,X(t,y))\ \hbox{ and }\ P(t,y):=Q(t,X(t,y)).
\end{equation}
We claim that System \eqref{NSE} recasts in 
\begin{equation}\label{NSL}
\begin{array}{lcr}
\eta_t=0 & \mbox{in} & \Omega \times (0,T), \\
\eta u_t - \nu \Delta_u u +\nabla_u P =0 & \mbox{in} & \Omega \times (0,T), \\
\divu u=0 & \mbox{in} & \Omega \times (0,T), \\
u=0 & \mbox{on} & \d \Omega \times (0,T), \\
u|_{t=t_0}=v|_{t=t_0}, \qquad \eta|_{t=0}=\rho|_{t=t_0}\  & \mbox{in} & \Omega,
\end{array}
\end{equation}
where operators $\Delta_u,\nabla_u, \divu$ correspond to the original operators $\Delta,\nabla,\div,$ respectively, 
after performing the change to the Lagrangian coordinates.  Index $u$ underlines the dependence with respect to $u$. 
Let us also notice that, as $v$ and $u$ vanish at the boundary, we do have  $X(t,\Omega)=\Omega$ for all $t$.

So let us now give a formal derivation of \eqref{NSL}.
First, given the definition of $X,$ it is obvious from the chain rule  that
$$
\d_t\eta(t,y)=(\d_t\rho+ v \cdot\nabla\rho)(t,x)
\quad\hbox{and}\quad
\d_t u(t,y)=(\d_t v+v\cdot\nabla v)(t,x)\quad\hbox{with }\ x:=X(t,y).
$$

The chain rule also yields
\begin{equation}
D_y P(t,y)=D_xQ(X(t,y))\cdot D_yX(t,y)
\quad\hbox{with }\  (D_yX)_{ij}:=\d_{y_j}X^i.
\end{equation}
Hence  we have 
\begin{equation}\label{eq:A}
D_x Q(t,x)=D_y P(t,y)\cdot A(t,y)\quad\hbox{with}\quad A(t,y):=(D_yX(t,y))^{-1}=D_xY(t,x).
\end{equation}

Next, we notice that if the transform $X$ is volume preserving 
then for any smooth enough vector-field $H$ we have
\begin{equation}\label{eq:div}
\divx  H(x)=\divy(A\bar H)(y)\quad\hbox{with }\ x=X(y)\ \hbox{ and }\ \bar H(y)=H(x).
\end{equation}
This stems from the following series of computations which uses 
the fact that $\det A\equiv1$ and the change of variable
$x=X(y)$: for any  smooth $q$ with compact support, we have
$$
\begin{array}{lll}
\Int q(x)\divx H(x)\,dx&=&-\Int D_xq(x)\cdot H(x)\,dx,\\[1ex]
&=&-\Int D_y\bar q(y)\cdot A(y)\cdot \bar H(y)\,dy,\\[1ex]
&=&\Int \bar q(y)\divy(A\bar H)(y)\,dy.\end{array}
$$

Combining \eqref{eq:A} and \eqref{eq:div}, we thus   deduce that, in Lagrangian coordinates 
operators $\nabla,$ $\div $ and $\Delta$ become
\begin{equation}
 \nabla_u := {}^T\!\!A \cdot \nabla_y, \quad \divu :=\div (A\cdot) \mbox{  ~~ and ~~ } \Delta_u:=\divu \nabla_u.
\end{equation}
In consequence, we have the following relations that will be of constant use:
\begin{eqnarray}\label{eq:duP}
&&(\nabla-\nabla_u)P=(\Id-{}^T\!A)\nabla P,\\\label{eq:Deltauu}
&&(\Delta-\Delta_u)u=\div((\Id-A{}^T\!A)\nabla u).
\end{eqnarray}

Let us finally emphasize that, owing to the chain rule, we   have
\begin{equation}\label{eq:magic}
\divy (A\cdot)=\divu=A:D_y.
\end{equation} 
This algebraic relation will be of fundamental importance in our analysis. 
\medbreak

The following statement ensures the  full equivalence  between \eqref{NSE} and \eqref{NSL}  under 
the assumptions of our results stated in Section \ref{s:main}.

\begin{prop}\label{p:change}
 Let  $1<p,q<\infty.$ Let   $\rho_0 \in L_\infty(\Omega)$ and 
$(u,P)$ be a solution to \eqref{NSL} such that   $u \in W^{2,1}_{q,p}(\Omega \times (0,T)),$  $\nabla P\in L_p(0,T;L_q(\Omega))$
and 
\begin{equation}
\int_0^T\|\nabla u\|_{L_\infty(\Omega)}\,dt\leq1/2.
\end{equation}
  Then 
$$
v(t,x)=u(t,y),\quad Q(t,x)=P(t,y) \mbox{ \ \ and \ \ } \rho(t,x)=\rho_0(y)
$$
with $x=X(t,y)$ given by \eqref{i3}
defines a $W^{2,1}_{q,p}$-solution to \eqref{NSE}.
\medbreak
Conversely, if  $\rho \in L_\infty(\Omega \times (0,T))$ and
 $(v,Q)$ with $v\in W^{2,1}_{q,p}(\Omega \times (0,T)),$    $\nabla v\in L_1(0,T;L_\infty(\Omega)),$ and 
 $\nabla Q\in L_p(0,T;L_q(\Omega))$
 is a solution to \eqref{NSE} then 
$$
u(t,y)=v(t,X(t,y)), \quad P(t,y)=Q(t,X(t,y)) \mbox{ \ \ and \ \  } \eta=\rho|_{t=0}
$$
defines a $W^{2,1}_{q,p}$-solution to \eqref{NSL}.
\end{prop}
\begin{p}
The proof goes along the lines of the corresponding one in the appendix of \cite{DM}.
Having $Du$ small enough in $L_1(0,T;L_\infty(\Omega))$ is of course fundamental. 
\end{p}


\section{Proof of Theorem \ref{th:uniq} -- uniqueness}\label{s:uniqueness} 

In this part we prove the uniqueness of solutions to System \eqref{NSL} under the assumptions of
Theorem \ref{th:uniq}. 
Here $\Omega$ is a $C^2$ bounded domain, or the whole space. 
 The proof is a
 straightforward application of Lemma \ref{l:l2stokes} to the equations in the Lagrangian form.
 The important fact is  that we have
$$\nabla v^i \in L_1(0,T;L_\infty(\Omega)) \qquad i=1,2.$$
Hence, taking $T$ small enough, one may assume with no loss of generality that 
\begin{equation}\label{eq:smallDu}\int_0^T\|\nabla v^i\|_{L_\infty(\Omega)}\,dt < \frac 12,\end{equation} so that
 Propositions \ref{p:lag} and \ref{p:change} apply. 
 In particular the regularity properties of 
those solutions in Lagrangian coordinates are the same as those  of Theorem \ref{th:uniq}. 
Hence it suffices to consider two solutions $(u^1,P^1)$ and $(u^2,P^2)$ to System \eqref{NSL} with the same initial data
and satisfying the conditions of Theorem \ref{th:uniq}.

Then, denoting $A^i:=A(u^i)$ (see \eqref{eq:A}),  $\du:=u^1-u^2$ and so on, we get 
\begin{equation}\label{NSL-uniq}
 \begin{array}{lcr}
  \eta\du_t - \nu \Delta\du + \nabla\dP= 
-\nu[(\Delta -\Delta_{u^1}) u^1 - (\Delta -\Delta_{u^2})u^2]\qquad&& \\
\qquad\qquad\qquad\qquad\qquad+
[(\nabla-\nabla_{u^1}) P^1 - (\nabla -\nabla_{u^2})P^2] & \mbox{in} & \Omega \times (0,T),\\[4pt]
\div\du=\div[( \Id -A^1)u^1 -(\Id -A^2) u^2] & \mbox{in} & \Omega \times (0,T),\\[4pt]
\du=0 & \mbox{ at } & \d \Omega \times (0,T),\\[4pt]
\du|_{t=0}=0 & \mbox{ in } & \Omega.
 \end{array}
\end{equation}
Let us underline that the boundary condition on $R$ from Lemma \ref{l:l2stokes} is fulfilled, since by definition $u^1$ and $u^2$ are zero at the boundary.
Therefore, keeping \eqref{eq:magic} in mind,  we obtain for some  constant $C$ depending  only on $\nu,$ $\inf \eta,$ $\sup\eta$ and  $\Omega$ 
 the inequality
\begin{equation}\label{uu8}
\|\du\|_{L_\infty(0,T;W^1_2(\Omega))}+ \|\du_t,\nabla^2\du,\nabla\dP\|_{L_2(\Omega \times (0,T))} \\\leq
C(I_1+I_2+I_3+I_4)
\end{equation}
with 
$$\begin{array}{lll}
I_1&:=& \| (\nabla-\nabla_{u^1})P^1 - (\nabla-\nabla_{u^2})P^2\|_{L_2(\Omega \times (0,T))}, \\[1ex]
I_2&:=&\| (\Id-A^1):Du^1 - (\Id -A^2):Du^2\|_{L_2(0,T;W^1_2(\Omega))},\\[1ex]
I_3&:=&\| [(\Delta -\Delta_{u^1})u^1 - (\Delta -\Delta_{u^2})u^2\|_{L_2(\Omega \times (0,T))}, \\[1ex]
I_4&:=&\| \d_t[(\Id-A^1)u^1 - (\Id-A^2)u^2]\|_{L_2(\Omega \times (0,T))}.
\end{array}
$$

In the following computations, we shall use repeatedly the fact that $\dA:=A^2-A^1$ satisfies 
\begin{equation}\label{eq:dA}
\dA(t)=\biggl(\int_0^tD\du\,d\tau\biggr)\cdot
\biggl(\sum_{k\geq1}\sum_{0\leq j<k} C_1^jC_2^{k-1-j}\biggr)
\quad\hbox{with}\quad
C_i(t):=\int_0^tDu^i\,d\tau.
\end{equation}

We concentrate on the case $n\geq3.$ We shall indicate how our arguments have to be modified
if $n=2,$ at the end of the section. 

In order to bound  $I_1,$ we write 
 \begin{equation}\label{u9}
 I_1(t)\leq 
\|(A^1-A^2) \nabla P^1\|_{L_2(\Omega \times (0,t))} +
\|(\Id - A^2) \nabla (P^1 - P^2)\|_{L_2(\Omega \times (0,t))}.
\end{equation}
It is clear that 
\begin{equation}
 \|(\Id - A^2) \nabla (P^1 - P^2)\|_{L_2(\Omega \times (0,t))} \leq Ct^{1/2} 
\|\nabla\dP\|_{L_2(\Omega \times (0,t))}\|Du^2\|_{L_2(0,t;L_\infty(\Omega))}.
\end{equation}
Let us notice that, according to  \eqref{i7},\eqref{eq:dA} and to the critical Sobolev embedding
of $W^1_2(\Omega)$ in $L_{2^*}(\Omega)$  (that is $1/2^*+1/n=1/2$), we have
$$\begin{array}{lll}
 \|\dA\|_{L_\infty(0,t;L_{2^*}(\Omega))} &\leq& C 
\|\Int_0^t |\nabla\du|dt'\|_{L_\infty(0,t;L_{2^*}(\Omega))}\\[1.5ex]
&\leq& Ct^{1/2}\|D^2 \du\|_{L_2(\Omega \times (0,t))}
\end{array}$$
with $C$ depending only on the norm of the two solutions on $[0,T].$
Therefore,  
$$ 
\begin{array}{lll}
 \|(A^1-A^2) \nabla P^1\|_{L_2(\Omega \times (0,t))}  &\leq& C 
 \|\dA\|_{L_\infty(0,t;L_{2^*}(\Omega))} 
\| \nabla P^1\|_{L_2(0,t;L_n(\Omega))}\\[1.5ex]
&\leq& C t^{1/2} \|D^2 \du\|_{L_2(\Omega\times(0,t))}\| \nabla P^1\|_{L_2(0,t;L_n(\Omega))}.
\end{array}
$$

Let us now  bound $I_2.$ Note that it suffices to bound the norm in 
$L_2(\Omega\times(0,T))$ of the gradient of the corresponding term. If $\Omega$ is
bounded this is a consequence of the Poincar\'e-Wirtinger inequality
as $\div\du$ has $0$ average over $\Omega,$ and if $\Omega=\R^n$ this stems from 
the fact that only the norm in $\dot W^1_2(\R^n)$ is involved (see Remark \ref{r:stokes}). 
Now, we notice   that
$$
(\Id-A^1):Du^1 - (\Id -A^2):Du^2=-\dA:D u^1+(A^2-\Id):D\du.
$$
First,  using the embedding of $W^1_2(\Omega)$ in $L_{2^*}(\Omega),$ 
 and keeping in mind \eqref{eq:smallDu} and that
\begin{equation}
Du^i\in L_2(0,T;L_\infty(\Omega))\ \hbox{ and }\ 
D^2u^i\in L_2(0,T;L_n(\Omega))\quad\hbox{for }\ i=1,2,
\end{equation}
we get for all $t\in[0,T],$
 $$
\begin{array}{lll}
\|D(\dA \!:\! Du^1)\|_{L_2(\Omega\times(0,t))} 
&\!\!\!\!\lesssim\!\!\!\!&\!\!\!\! \Bigl\| |Du^1|\!\Int_0^\tau \!|D^2\du|\,d\tau' \Bigr\|_{L_2(\Omega\times (0,t))}
 \!+\! \Bigl\||D^2u^1|\!\int_0^\tau \!|D\du|\,d\tau' \Bigr\|_{L_2(\Omega\times (0,t))}
\\[1.5ex]
&\!\!\!\!\lesssim\!\!\!\!& t^{1/2}\Bigl(\|D^2\du\|_{L_2(\Omega\times(0,t))} \|Du^1\|_{L_2(0,t; L_\infty(\Omega))} 
 \\&&\qquad\qquad+ \|D\du\|_{L_2(0,t;L_{2^*}(\Omega))} \|D^2u^1\|_{L_2(0,t;L_n(\Omega))}\Bigr) 
\\[1ex] &\!\!\!\!\lesssim\!\!\!\! &t^{1/2} \|D^2\du\|_{L_2(\Omega\times(0,t))}.
\end{array}
$$
Second, we have
$$
\begin{array}{lll}
\|D((A^2-\Id)\!:\!D\du)\|_{L_2(\Omega\times(0,t))} 
&\!\!\!\!\lesssim\!\!\!\!&\!\! 
\|D A^2\!\otimes\!  D\du\|_{L_2(\Omega\times(0,T))}+ \|(A^2-\Id)\!\otimes\! D^2\du)\|_{L_2(\Omega\times(0,t))} 
\\[1ex]
&\!\!\!\!\lesssim\!\!\!\!& t^{1/2}\Bigl(\|D^2u^2\|_{L_2(0,t;L_n(\Omega))} 
\|D\du\|_{L_2(0,t; L_{2^*}(\Omega))} \\&&+\|Du^2\|_{L_1(0,t;L_\infty(\Omega))} \|D^2\du\|_{L_2(\Omega\times(0,t))}\Bigr).
\end{array}
$$
So finally for all $t\in[0,T],$
\begin{equation}\label{eq:I4}
I_2(t)\leq Ct^{1/2} \|\nabla^2\du\|_{L_2(\Omega\times(0,t))}
\end{equation}
with $C$ depending only the norm of the solutions over $[0,T].$
The term $I_3$ may be handled 
along the same lines.  Indeed we have
$$
I_3(t)=\Bigl\|\div\Bigl(\bigl(\Id-A^1\,{}^T\!A^1\bigr)\nabla u^1-\bigl(\Id-A^2\,{}^T\!A^2\bigr)\nabla u^2\Bigr)\Bigr\|_{L_2(\Omega\times(0,t))}.
$$
Finally, we examine $I_4.$ 
Using again \eqref{eq:dA}, we get (with the convention that $Du^{1,2}$ denotes the components
of $Du^1$ and $Du^2$):
$$\begin{array}{lll}
\!\!\|\d_t[\dA \,u^1]\|_{L_2(\Omega\times(0,t))} 
&\!\!\!\!\lesssim\!\!\!& \|D\du\, u^1\|_{L_2(\Omega\times(0,t))}
+\Bigl\|\Int_0^\tau |D\du|\,d\tau'|Du^{1,2}|\,|u^1|\Bigr\|_{L_2(\Omega\times(0,t))}\\[1.5ex] 
&&\qquad\qquad\qquad\qquad+ \Bigl\|\Int_0^\tau |D\du|\,d\tau' \, |u^1_t| \Bigr\|_{L_2(\Omega\times(0,t))}\\[3ex]
&\!\!\!\!\lesssim\!\!\!& \|D\du\|_{L_\infty(0,t;L_2(\Omega))} \|u^1\|_{L_2(0,t;L_\infty(\Omega))}

 \\&&\!\!\!\!\!\!\!\!\!\qquad\!+t^{1/2}\|D\du\|_{L_2(0,t;L_{2^*}(\Omega))}
\bigl(\|u_t^1\|_{L_2(0,t;L_n(\Omega))}\!+\!\||u^1||Du^{1,2}|\|_{L_2(0,t;L_n(\Omega))}\bigr)
\\[2ex]
&\!\!\!\!\lesssim\!\!\!& t^{1/2} \|D^2\du\|_{L_2(\Omega\times(0,t))}+\ep(t) \|D\du\|_{L_\infty(0,t;L_2(\Omega))}, 
\end{array}
$$
with $\lim_{t\rightarrow0}\ep(t)=0$ because
\begin{equation}\label{eq:u8}
\d_tu^1,\  u^1\otimes Du^1\hbox{ and } u^2\otimes Du^1\ \hbox{ are in }\ 
 L_2(0,T;L_n(\Omega)).
\end{equation}
At the same time, we have, for all $t\in[0,T],$
$$
\begin{array}{lll}
\|\d_t((\Id-A^2)\du)\|_{L_2(\Omega\times(0,t))}&\!\!\!\lesssim\!\!\!&\|Du^2\du\|_{L_2(\Omega\times(0,t))}+\|(\Id-A^2)\d_t\du\|_{L_2(\Omega\times(0,t))}\\[1ex]
&\!\!\!\lesssim\!\!\!&\|\du\|_{L_\infty(0,t;L_{n^*}(\Omega))}
\|Du^2\|_{L_2(0,t;L_n(\Omega))}\\&&\qquad\qquad
+\|Du^2\|_{L_1(0,t;L_\infty(\Omega))}\|\d_t\du\|_{L_2(\Omega\times(0,t))}.
\end{array}
$$
So one may conclude that  
$$
I_4(t)\leq t^{1/2} \|D^2\du\|_{L_2(\Omega\times(0,t))}
+\ep(t)\bigl( \|D\du\|_{L_\infty(0,t;L_2(\Omega))}+\|\d_t\du\|_{L_2(\Omega\times(0,t))}\bigr). 
$$

So finally in the  case  $n\geq 3$, putting together all the previous
inequalities yields  for all $t\in (0,T),$
$$
 \displaylines{\|\du\|_{L_\infty(0,t;W^1_2(\Omega))}+\|\du_t,\nabla^2\du,\nabla\dP\|_{L_2(\Omega \times (0,t))} \\
\hfill\cr\hfill\leq \ep(t)\Bigl(\|\du\|_{L_\infty(0,t;W^1_2(\Omega))}+\|\du_t,\nabla^2\du,\nabla\dP\|_{L_2(\Omega \times (0,t))}\Bigr)}
$$
for some positive function $\ep$ going to $0$ at $0.$  
Uniqueness follows on a sufficiently small time interval, then on the whole interval $[0,T]$
thanks to a standard connectivity (or bootstrap) argument.
\medbreak
Let us now explain how the arguments have to be modified 
in the two-dimensional case.  One cannot follow exactly 
the above approach owing to the failure of the embedding of $W^1_2(\Omega)$ 
in  $L_\infty(\Omega)$. So we have
to assume slightly higher regularity, namely $\nabla P^1,\nabla P^2 \in L_2(0,T;L_q(\Omega))$ with $q>2,$
and so on. For instance, setting  $m\in(2,\infty)$ such 
 that $1/m+1/q=1/2,$  we may write
$$
\begin{array}{lll}
 \|(A^1-A^2)\nabla P^1\|_{L_2(\Omega \times (0,t))} &\leq& Ct^{1/2} \|\nabla\du\|_{L_2(0,t;L_m(\Omega))} \|\nabla P^2\|_{L_2(0,t;
L_q(\Omega))}\\[1.5ex]
&\leq& C t^{1/2} \|D^2 \du\|_{L_2(\Omega\times(0,t))} \|\nabla P^2\|_{L_2(0,t;
L_q(\Omega))}.
\end{array}
$$
The other terms of \eqref{uu8} may be handled similarly. The details are left to the reader.
 Theorem \ref{th:uniq} is thus proved. 

\begin{rem}
 Here we would like to explain the reason why we use the $W^{2,1}_{2,2}$ regularity for the velocity
 to establish uniqueness. 
Concentrate our attention on $n=3$. A direct $L_2$-energy method (i.e. 
 testing \eqref{NSL-uniq} by $\du$)  requires our bounding 
$(\nabla - \nabla_{u^1})P^1-(\nabla-\nabla_{u^2})P^2$ in $L_1(0,T;L_2(\Omega))$,
hence the following computation:
$$
\begin{array}{lll}
\left| \Int_0^T \int_\Omega \dA\nabla P^1 \,\du \,dx\, dt \right|  
&\leq& 
C\|\dA\|_{L_\infty(0,T;L_2(\Omega))} \|\nabla P^1\|_{L_2(0,T;L_3(\Omega))} \|\du\|_{L_2(0,T;L_6(\Omega))}
\\ &\leq& 
CT^{1/2}   \|\nabla P^1\|_{L_2(0,T;L_3(\Omega))} \|\nabla \du \|_{L_2(\Omega\times(0,T))}.
\end{array}
$$
So we need  $\nabla P^1 \in L_2(0,T;L_3(\Omega))$ which is naturally related to $u^1 \in W^{2,1}_{3,2}$. In addition  integrating  by parts 
 in the left-hand side  of the above inequality, we need to keep track of
 $\nabla^2 \du$ as well as of $\nabla\dP$ in $L_2(\Omega \times (0,T)).$
 Those two terms are  out of  control if resorting only to the basic energy inequality.
\end{rem}


\section{Proof of Theorem \ref{th:exist} -- existence} \label{s:existence}

The uniqueness property of the system is important, but to have the full picture
 of the well-posedness issue, we now have to show that there exist solutions \emph{with merely bounded  density}
 for which Theorem~\ref{th:uniq} applies. With the method that is proposed below, 
 much more regularity is needed for the initial velocity. However the assumption over the initial density stays that same:  it just has to be bounded and bounded away from zero.

\subsection{A priori estimates}

We first  concentrate on the proof of a priori estimates for 
a smooth solution $(u,P)$ to  \eqref{NSL}. 
To simplify the presentation,  we consider the case where $\Omega$ is a $C^2$ bounded
domain of $\R^n.$ The whole space case may be achieved by similar arguments : this is just
a matter of using homogeneous norms $\|\cdot\|_{\dot W^1_2(\R^n)}$ and 
$\|\cdot\|_{\dot B^{2-2/n^*}_{n^*,n^*}(\R^n)}$  and resorting to Remark \ref{r:stokes}. 
\smallbreak
In order  to prove a priori estimates for $(u,P),$ let us assume in addition that
$T$ has been chosen so that  (say)
\begin{equation}\label{eq:smallu}
\int_0^T\|\nabla u\|_{L_\infty(\Omega)}\,dt\leq1/2.
\end{equation}
This enables us to go from \eqref{NSE} to \eqref{NSL} (and conversely).
For any (possibly large) initial velocity $v_0\in B^{2-2/n^*}_{n^*,n^*}(\Omega),$
and $\rho_0\in L_\infty(\Omega)$ bounded away from zero, 
 we want to find a bound  for a solution $(u,P)$ given by Theorem \ref{th:exist}. 
In other words, we want to control the following quantity:
\begin{equation}
 \Xi_{(u,P)}(T):=  \|u_t\|_{L_\infty(0,T;L_2(\Omega))} + \|\nabla u_t\|_{L_2(\Omega \times (0,T))} +
\|u_t,\nabla^2u,\nabla P\|_{L_{n^*}(\Omega \times (0,T))},
\end{equation}
with $n^*=2\bigl(\frac{n+2}{n}\bigr),$ if  $T$ is small enough.

\medbreak
Let us first notice that, by standard Sobolev embedding
\begin{equation}\label{eq:smallu1}
\int_0^T\|\nabla u\|_{L_\infty(\Omega)}\,dt\leq CT^{1-\frac{1}{n^*}} \Xi_{(u,P)}(T).
\end{equation}
which guarantees \eqref{eq:smallu} for small times.

\medbreak
In order to use Lemma \ref{l:l2stokes} we restate System \eqref{NSL} as follows
(of course $\eta=\rho_0$ and $u_0=v_0$):
\begin{equation}\label{NSL-1}
\begin{array}{lcr}
\eta u_t - \nu \Delta u +\nabla P =-\nu(\Delta -\Delta_u)u+(\nabla-\nabla_u)P \qquad & \mbox{in} & \Omega \times (0,T), \\[4pt]
\div u=\div\bigl((\Id-A)u\bigr) & \mbox{in} & \Omega \times (0,T), \\[4pt]
u=0 & \mbox{on} & \d \Omega \times (0,T), \\[4pt]
u|_{t=0}=u_0 & \mbox{in} & \Omega.
\end{array}
\end{equation}
Then keeping \eqref{eq:duP}, \eqref{eq:Deltauu}, \eqref{eq:magic}, and Proposition \ref{p:lag} in mind, we get
 for some constant $C=C(\nu,\Omega),$
\begin{multline}\label{p1}
\sup_{0\leq t \leq T}\|u(t)\|_{W^1_2(\Omega)}+ \|u_t,\nabla^2u,\nabla P \|_{L_2(\Omega\times(0,T))} \\
\leq C\big( \|\Id-A\|_{L_\infty(\Omega \times (0,T))} \|u_t,\nabla^2u,\nabla P \|_{L_2(\Omega \times (0,T))} \\
+ \|\nabla A \, \nabla u, A_t\,  u \|_{L_2(\Omega\times(0,T))} + \|u_0\|_{W^1_2(\Omega)}\big).
\end{multline}
The $W^{2,1}_{2,2}(\Omega \times (0,T))$ regularity of the velocity, coming from \eqref{p1}, is not sufficient to control the Lagrangian coordinates, namely the terms containing $A$ in the right-hand side of \eqref{p1},
because  $\nabla W^{2,1}_{2,2}(\Omega \times (0,T))$ \emph{is not} embedded in $L_1(0,T;L_\infty(\Omega)).$
 Hence, to close the estimates,  higher regularity is needed.
Differentiating \eqref{NSL-1} once with respect to time is the easiest way to achieve it, because it does not affect
 the irregular density which is time independent in the Lagrangian setting. We get
\begin{equation}\label{NSL-2}
\begin{array}{lcr}
\eta u_{tt} - \nu \Delta u_t +\nabla_u P_t = && \\
\qquad\qquad -\nu(\Delta -\Delta_u)u_t+ \nu(\Delta_u)_tu-(\nabla_u)_t P \qquad \qquad & \mbox{in} & \Omega \times (0,T), \\[4pt]
\divu u_t= -\div A_t u & \mbox{in} & \Omega \times (0,T), \\[4pt]
u_t=0 & \mbox{on} & \d \Omega \times (0,T).
\end{array}
\end{equation}

 At this stage  the question of the regularity of $u_t|_{t=0}$ arises. 
This information can be found out  only from the equations. At time  $t=0$ the Eulerian and Lagrangian coordinates coincide (that  is $A=\Id$), 
so the regularity of $u_t|_{t=0}$ is just that of $\eta^{-1}(\nabla P-\nu\Delta u)|_{t=0}.$
However the regularity of $\nabla  P|_{t=0}$ is unknown, so we rather have to use the fact that
 differentiating  $\eqref{NSL-1}_2$ with respect to  $t$ implies that 
\begin{equation}\label{in1}
 \eta u_t|_{t=0} + \nabla P|_{t=0}= \nu \Delta u_0, \qquad \qquad \div u_t|_{t=0} = -\div (A_t|_{t=0} u_0).
\end{equation}
Note that
$A_t|_{t=0}$ need not be trivial so  in order to bound  $u_t|_{t=0}$ in $L_2(\Omega),$
we first have to remove its potential part.  For that, we 
use  the Bogovski\u{\i} operator $\cB$ (see Lemma \ref{l:bogovskii}) setting 
$$
 \phi = {\cal B}[-\div (A_t|_{t=0} u_0)] \mbox{ ~~~ so that ~~~ } \div \phi = - \div (A_t|_{t=0} u_0) \mbox{ in } \Omega, \quad 
\phi =0 \mbox{ at } \d \Omega.
$$
Let us notice that, because 
$$
A(t,y)\cdot DX(t,y)=\Id\ \hbox{ and }\ DX|_{t=0}=\Id\quad\hbox{with}\quad
X(t,y)=\Id+\int_0^t u(\tau,y)\,d\tau,
$$
we have  $A_t|_{t=0}=-Du_0,$
hence  $A_t|_{t=0} u_0 = -(u_0\cdot\nabla u_0)=-\div(u_0\otimes u_0).$
\smallbreak
Now, $W^2_2(\Omega)$ is an algebra if $n=2,3.$ Hence $A_t|_{t=0} u_0$ is in $W^1_2(\Omega)$
and the function $\phi$ defined above
 is in $W^1_2(\Omega)$\footnote{In fact, the function $\phi$ is  in $W^2_2(\Omega)$
but we shall not take advantage of this in what follows.} and satisfies:
\begin{equation}\label{in4}
 \|\phi \|_{W^1_2(\Omega)} \leq C\|u_0\|^2_{W^2_2(\Omega)}.
\end{equation}
Therefore System \eqref{in1} recasts in
$$\begin{array}{l}
 \eta (u_t|_{t=0}-\phi) + \nabla P|_{t=0}= \nu \Delta u_0-\eta \phi\  \mbox{ in }\ \Omega\\[1ex]
  \div (u_t|_{t=0} -\phi)=0\  \mbox{ in }\ \Omega\\[1ex]
(u_t|_{t=0} -\phi)|_{\d\Omega}=0\  \mbox{ on }\ \d\Omega.
\end{array}
$$
Now, testing  the first equation by $(u_t|_{t=0} - \phi)$ we get:
\begin{equation}\label{in6}
 \int_\Omega \eta \bigl|u_t|_{t=0} -\phi\bigr|^2\, dx \leq
 \int_\Omega \eta^{-1}\bigl|\nu\Delta u_0-\eta\phi\bigr|^2\,dx.
\end{equation}
Thus, due to \eqref{in4}, we discover that $u_t|_{t=0}$ is in $L_2(\Omega)$ and that
\begin{equation}\label{in7}
 \|u_t|_{t=0} \|_{L_2(\Omega)} \leq C_{\nu,\eta}(\|u_0\|_{W^2_2(\Omega)} + \|u_0\|^2_{W^2_2(\Omega)}).
\end{equation}

At this point, we would like to apply  an energy method to \eqref{NSL-2}. However, as $\divu u_t$ may be nonzero,
one cannot
eliminate the term coming from $P_t$ (which is out of control). 
So we modify $\eqref{NSL-2}_2$  by  introducing a vector-field $\xi$ so that 
\begin{equation}\label{in8}
 \begin{array}{lr}
 \divu \xi = -\div (A_t u) \;\;  &\hbox{in }\ \Omega, \\
\xi =0 &\hbox{on }\  \d \Omega.
\end{array}
\end{equation}
We need $\xi$ to satisfy suitable estimates (in terms of the right-hand side) 
in $L_\infty(0,T;L_2(\Omega))\cap L_2(0,T;W^1_2(\Omega))$
and  $\xi_t$ to be bounded in  $L_2(\Omega\times(0,T)).$ 
This may be done by  means of a Bogovski\u{\i} type operator construction as in \cite{DM-luminy}. 
Here we shall define $\xi$ (treating $t$ as parameter) according to   Lemma \ref{l:divA} in the Appendix.
\smallbreak
 Let us start with the bound in $L_\infty(0,T;L_2(\Omega))$: 
we have
$$\|\xi\|_{L_\infty(0,T;L_2(\Omega))}\lesssim \|A_tu\|_{L_\infty(0,T;L_2(\Omega))}.
$$
Therefore, using the fact that
\begin{equation}\label{eq:At}
A_t=\biggl(\sum_{k\geq0}(k+1)(-1)^{k+1}\biggl(\int_0^tD_yu\,d\tau\biggr)^k\biggr)\cdot D_yu,
\end{equation}
we get (remember \eqref{eq:smallu})
\begin{eqnarray}\label{in8a}
 &&\|\xi\|_{L_\infty(0,T;L_2(\Omega))}\lesssim\| u \otimes \nabla u  \|_{L_\infty(0,T;L_2(\Omega))}\nonumber \\
&&\phantom{\|\xi\|_{L_\infty(0,T;L_2(\Omega))}}\lesssim\|u\|_{L_\infty(0,T;L_\infty(\Omega))} \|\nabla u \|_{L_\infty(0,T;L_2(\Omega))}. 
\end{eqnarray}

In order to bound  the right-hand side of \eqref{in8a}, we apply 
the following classical parabolic estimate (which is related to our definition of Besov
spaces in \eqref{xx1}):
\begin{equation}\label{eq:parabolic}
 \|u\|_{L_\infty(0,T;B^{2-2/p}_{p,p}(\Omega))} \leq C( \|u_0\|_{B^{2-2/p}_{p,p}(\Omega)} 
+\|u_t,\nabla^2 u\|_{L_p(\Omega \times (0,T))}).
\end{equation}
Now, owing to Sobolev embedding, it is clear that the left-hand side  of \eqref{eq:parabolic} controls the $L_\infty$ norm 
whenever $2-2/p> n/p,$ that is  $p>(n+2)/2.$ The constant in \eqref{eq:parabolic} is time independent. 
Therefore, for any $m\in((n+2)/2,n^*),$ we have
\begin{eqnarray}\label{in9}
&&\|u\|_{L_\infty(\Omega \times (0,T))}  \leq C_m\bigl(
\|u_0\|_{B^{2-2/m}_{m,m}(\Omega)} + \|u_t,\nabla^2 u\|_{L_{m}(\Omega \times (0,T))}\bigr)\nonumber
\\ &&\phantom{\|u\|_{L_\infty(\Omega \times (0,T))}}\leq C_m\bigl(\|u_0\|_{B^{2-2/{n^*}}_{n^*,n^*}(\Omega)} + T^{\frac 1m - \frac{1}{n^*}} \,\Xi_{(u,P)}(T)\bigr).
\end{eqnarray}
Inequality \eqref{eq:parabolic} with $p=2$ also yields
\begin{eqnarray}\label{in10}
&& \|\nabla u\|_{L_\infty(0,T;L_2(\Omega)} \lesssim
 \|u_0\|_{W^1_2(\Omega)} + \|u_t,\nabla^2 u\|_{L_2(\Omega \times (0,T))}\nonumber\\
 &&\phantom{ \|\nabla u\|_{L_\infty(0,T;L_2(\Omega)} }\lesssim 
\|u_0\|_{W^1_2(\Omega)} + T^{\frac 12 - \frac{1}{n^*}} \Xi_{(u,P)}(T).
\end{eqnarray}
 So, putting  \eqref{in8a}, \eqref{in9} and \eqref{in10} together, we get  for some $\delta>0,$
 \begin{equation}\label{in10a}
 \|\xi\|_{L_\infty(0,T;L_2(\Omega))}\lesssim \|u_0\|^2_{W^2_2(\Omega)} +  T^\delta \, \Xi_{(u,P)}^2(T).
\end{equation}

Next, in order to  bound $\xi$ in $L_2(0,T;W^1_2(\Omega)),$ we use  the fact that 
$$\divu\xi=-A_t:Du.$$
In effect, owing to \eqref{eq:div}, one may write
$\div(Au_t)=A:Du_t,$
hence the above relation may be obtained by taking the time derivative of
$\div(Au)=A:Du.$
So,  using  Lemma \ref{l:divA} and 
remembering that $n^* >n,$ we get, for some $\delta>0,$
$$\begin{array}{lll}
 \|\xi \|_{L_2(0,T;W^1_2(\Omega))} &\!\!\!\leq\!\!\!& C\| |\nabla u |^2 \|_{L_2(\Omega\times(0,T))}\\[1ex]
 &\!\!\!\leq\!\!\!& C \|\nabla u\|_{L_\infty(0,T;L_2(\Omega))}\|\nabla u\|_{L_2(0,T;L_\infty(\Omega))}\\[1ex]
&\!\!\!\leq\!\!\!& C(\|\nabla u_0\|_{L_2(\Omega)} + T^{1/2}\|\nabla u_t\|_{L_2(\Omega\times(0,T))}) 
\|\nabla u\|_{L_2(0,T; W^1_{n^*}(\Omega))}\\[1ex]
&\!\!\!\leq\!\!\!& C T^{\delta}\Xi_{(u,P)}(T)(\|u_0\|_{W^1_2(\Omega)} + T^{\delta}\Xi_{(u,P)}(T)).
\end{array}
$$
Finally, let us   bound $\xi_t$ in  $L_2(\Omega\times(0,T)).$
For that, we  apply the last part of Lemma \ref{l:divA} which yields
\begin{equation}\label{in11}
 \|\xi_t \|_{L_2(\Omega\times(0,T))} \leq C \|A_t \xi, A_{tt} u, A_t u_t\|_{L_2(\Omega\times(0,T))}.
\end{equation}
Using \eqref{in8a}, \eqref{in9} and \eqref{in10a}, we get
$$
\begin{array}{lll}
 \|A_t \xi\|_{L_2(\Omega\times(0,T))} &\leq &C\|\nabla u\|_{L_2(0,T;L_\infty(\Omega))} \|\xi\|_{L_\infty(0,T;L_2(\Omega))}\\[1ex]&\leq&
C\bigl(\|u_0\|_{W^{2}_{2}(\Omega)}^2 +T^\delta \Xi^2_{(u,P)}(T)\bigr)
\|\nabla u\|_{L_2(0,T;L_\infty(\Omega))},
\end{array}
$$
$$
\begin{array}{lll}
 \|A_{tt} u \|_{L_2(\Omega\times(0,T))}&\leq& \|\nabla  u_t\|_{L_2(\Omega\times(0,T))} \|u\|_{L_\infty(0,T;L_\infty(\Omega))}\\[1ex]
&\leq& C\|\nabla u_t\|_{L_2(\Omega\times(0,T))} (\|u_0\|_{W^{2-2/n^*}_{n^*}(\Omega)} +T^\delta \Xi_{(u,P)}(T)),
\end{array}
$$
$$\displaylines{\quad
 \|A_t u_t\|_{L_2(\Omega\times(0,T))} \leq C\|\nabla u\otimes  u_t\|_{L_2(\Omega\times(0,T))} \hfill\cr\hfill
\leq 
C\|\nabla u \|_{L_2(0,T;L_\infty(\Omega))} (\|u_t-\xi\|_{L_\infty(0,T;L_2(\Omega))} +
\|\xi\|_{L_\infty(0,T;L_2(\Omega))}).\quad}
$$
So one may conclude that
\begin{multline}\label{in15}
\|\xi_t\|_{L_2(\Omega\times(0,T))}\lesssim (\|u_0\|_{W^2_2(\Omega)}^2 + T^{\delta}\Xi_{(u,P)}^2(T))
\|\nabla u\|_{L_2(0,T;L_\infty(\Omega))} +\|u_0\|_{W^{2-2/n^*}_{n^*}(\Omega)}\\ + T^{\delta}\Xi_{(u,P)}(T))\|\nabla u_t\|_{L_2(0,T;L_2(\Omega))}+
\|\nabla u\|_{L_2(0,T;L_\infty(\Omega))} \|u_t-\xi\|_{L_\infty(0,T;L_2(\Omega))}.
\end{multline}
Note that in \eqref{in15},  there is no factor $T^\delta$ for the leading order
terms   $\|u_t-\xi\|_{L_\infty(0,T;L_2(\Omega))}$ and  $\|\nabla u_t\|_{L_2(\Omega\times(0,T))}.$ 
\medbreak
Once the vector-field $\xi$  has been constructed, one may  recast System \eqref{NSL-2} in
\begin{equation}\label{NSL-2a}
\begin{array}{lcr}
\eta (u_{t}-\xi)_t - \nu \Delta (u_t-\xi) +\nabla_u P_t \qquad&& \\
\qquad=-\nu(\Delta -\Delta_u)u_t+ \nu( \Delta_u)_tu-(\nabla_u)_t P
-\eta \xi_t + \nu \Delta \xi & \mbox{in} & \Omega \times (0,T), \\[4pt]
\divu (u_t-\xi)= 0 & \mbox{in} & \Omega \times (0,T), \\[4pt]
u_t-\xi =0 & \mbox{on} & \d \Omega \times (0,T), \\[4pt]
(u_t-\xi)|_{t=0} \in L_2(\Omega)
& \mbox{in} & \Omega.
\end{array}
\end{equation}

Note that, now, $\divu(u_t-\xi)=0$ and that $(u_t-\xi)|_{t=0}$ is in $L_2(\Omega)$ with
\begin{equation}\label{in13}
\|(u_t-\xi)|_{t=0}\|_{L_2(\Omega)}\leq C\bigl(\|u_0\|_{W^2_2(\Omega)}+\|u_0\|_{W^2_2(\Omega)}^2\bigr).
\end{equation}
  So taking the $L_2(\Omega)$-inner product of $(\ref{NSL-2a})_1$ with  $u_t-\xi,$  there is no term
generated by $\nabla_uP_t$ and we thus get
$$
 \frac 12 \int_\Omega \eta |u_t -\xi|^2 dx\Big|_{t=T} + \nu \int_0^T\!\!\!\int_\Omega |\nabla (u_t -\xi)|^2 dxdt \leq
 \frac 12 \int_\Omega \eta  |u_t -\xi|^2\,dx\Big|_{t=0}+\sum_{j=1}^5I_j
 $$
 with 
 $$
 \begin{array}{lll}
 I_1&:=& C\nu\Int_0^T\!\!\! \int_\Omega|(A -\Id)|\,|\nabla u_t |\, |\nabla (u_t-\xi) |\,dx\,dt,\\[2ex]
  I_2&:=&C \nu\Int_0^T \!\!\!\int_\Omega  |\nabla u|^2| \nabla (u_t-\xi)| \,dx\,dt,\\[2ex]
 I_3&:=& C\Int_0^T \!\!\!\int_\Omega |\nabla u | \,|\nabla P |\, |u_t-\xi|\,dx\,dt,\\[2ex]
  I_4&:=& \Int_0^T \!\!\!\int_\Omega \eta|\xi_t \cdot (u_t-\xi)|\,dx\,dt,\\[2ex]
  I_5&:=& \nu\Int_0^T \!\!\!\int_\Omega  |\nabla \xi \cdot \nabla (u_t -\xi)| \, dx\,dt.
  \end{array}
  $$

In order to bound terms $I_1,$ $I_2,$ $I_3$  and  $I_5,$ we use
H\"older and Young inequalities. We get for all $\ep>0,$ 
\begin{eqnarray}\label{in17}
 &&I_1 \leq \ep\nu \|\nabla(u_t-\xi)\|_{L_2(\Omega \times(0,T))}^2
+C_{\ep,\nu}\|A-\Id\|_{L_\infty(\Omega \times(0,T))}^2 \|\nabla u_t\|_{L_2(\Omega\times(0,T))}^2,\quad \\
\label{in18}
&&I_2 \leq \ep \nu \|\nabla(u_t-\xi)\|_{L_2( \Omega \times (0,T))}^2
+C_{\ep,\nu} \| |\nabla u|^2 \|_{L_2(\Omega\times(0,T))}^2,\\\label{in19}
 &&I_3\leq \ep\|u_t-\xi\|_{L_\infty(0,T;L_2(\Omega))}^2 +C_\ep\|\nabla u\|_{L_2(0,T;L_\infty(\Omega))}^2 \|\nabla P\|^2_{L_2(\Omega\times(0,T))},\\
&&I_5\leq \ep \nu \|\nabla (u_t -\xi)\|_{L_2(\Omega \times (0,T))}^2+C_{\ep,\nu} \|\nabla \xi\|_{L_2(\Omega \times(0,T))}^2.
\end{eqnarray}

Inequality \eqref{in19} deserves a remark : in order to ``close the estimates'',  we have to factor out the last term in the right-hand side by a quantity which is small enough
when $T$ goes to $0.$  Here this follows from the embedding 
$W^{2,1}_{n^*,n^*} \subset L_2(0,T;L_\infty(\Omega))$ which gives,
because $n^*>2,$ 
\begin{equation}
 \|\nabla u\|_{L_2(0,T;L_\infty(\Omega))} \leq CT^{ \frac 12 - \frac{1}{n^*} } \|u\|_{W^{2,1}_{n^*,n^*}(\Omega \times (0,T))}.
\end{equation}
Finally, taking $m\in(1,2)$ so that  $1=\frac 1m + \frac{1}{n^*},$ and $\delta := \frac{2}{m} -1,$ 
we may write
$$
\begin{array}{lll}
I_4&\leq&  \|u_t-\xi\|_{L_{n^*}(\Omega \times (0,T))}\|\xi_t\|_{L_{m}(\Omega \times (0,T))},\\[1ex]
& \leq&  \ep \|u_t-\xi\|_{L_{n^*}(\Omega \times (0,T))}^2 +C_\ep T^\delta \|\xi_t\|_{L_{2}(\Omega \times (0,T))}^2.
\end{array}
$$
Combining interpolation and Sobolev embedding, we may write for all $p\in(2,\infty),$
$$
\|u_t-\xi\|_{L_p(0,T;L_q(\Omega))} 
\lesssim \|u_t-\xi\|_{L_\infty(0,T;L_2(\Omega))}^{1-2/p}\|D(u_t-\xi)\|_{L_2(\Omega\times(0,T))}^{2/p}, 
$$
with $n/q=n/2-2/p.$ So taking $p=q=n^*:=2(n+2)/n,$ we get
\begin{equation}\label{in20a}
 \|u_t-\xi\|_{L_{n^*}(\Omega \times (0,T))}\leq C(\|u_t-\xi\|_{L_\infty(0,T;L_2(\Omega))}
+\|\nabla (u_t-\xi)\|_{L_2(\Omega \times (0,T))}).
\end{equation}
Therefore, the above estimates for $I_1$ to $I_5$ (with $\ep$ small enough) eventually imply that 
$$
\displaylines{
\|u_t-\xi\|_{L_\infty(0,T;L_2(\Omega))}^2+\nu\|\nabla(u_t-\xi)\|_{L_2(\Omega\times(0,T))}^2
\leq 2\|(u_t-\xi)|_{t=0}\|_{L_2(\Omega)}^2\hfill\cr\hfill
+C\Bigl(\|\nabla u\|_{L_2(0,T;L_\infty(\Omega))}^2\bigl(\|\nabla u\|_{L_\infty(0,T;L_2(\Omega))}^2
+\|\nabla P\|_{L_2(\Omega \times (0,T))}^2\bigr)
\hfill\cr\hfill+\|\nabla u\|_{L_1(0,T;L_\infty(\Omega))}^2\|\nabla u_t\|_{L_2(\Omega \times (0,T))}^2
+T^\delta \|\xi_t\|_{L_2(\Omega \times (0,T))}^2+\|\nabla\xi\|^2_{L_2(\Omega \times (0,T))}\Bigr),}
$$
whence, using also the estimates for $\xi,$ $\xi_t$ and for $\nabla u$ in $L_2(0,T;L_\infty(\Omega)),$
we end up with \begin{equation}\label{in20}
 \|u_t\|_{L_\infty(0,T;L_2(\Omega))}+  \|\nabla u_t \|_{L_2(\Omega \times (0,T))}
\leq C_{u_0}\bigl(1+T^\delta \Xi_{(u,P)}(T)\bigr)^3.
\end{equation}
Here $\delta>0.$ Let us also stress that   $C_{u_0}$ depends only on $\|u_0\|_{W^2_2(\Omega)},$
$\|\rho_0\|_{L_\infty(\Omega)},$ $\|\rho_0^{-1}\|_{L_\infty(\Omega)},$ 
$\Omega$ and $\nu.$ In particular, it is time-independent.

\begin{rem}
At  this stage,  we find a limitation on the  dimension of the domain: 
as we need to have $\nabla u\in L_{2}(0,T;L_\infty(\Omega))$,  embedding requires that  $n^*>n$. 
This is fulfilled if $n=2$ (because $2^*=4$) or $n=3$ (because $3^*=10/3$) but this is no
longer satisfied in  higher dimension. So we see that our method cannot 
be directly applied for $n \geq 4$ unless we differentiate the system with respect to time, more times. 
Physical motivation for considering dimension $n\geq4$ is
unclear, though.
\end{rem}

Keeping \eqref{in20a} in mind, we see that in order to close the estimates, 
it suffices  to bound the terms $\|\nabla^2u\|_{L_{n^*}(\Omega \times (0,T))}$
and $\| \nabla P\|_{L_{n^*}(\Omega \times (0,T))}$ which appear in the right-hand side of \eqref{in20}.
For that, we rewrite System \eqref{NSL-1} as a \emph{stationary} Stokes system, treating  $\eta u_t$ as a source term, 
and the time variable as a parameter. So we consider 
\begin{equation}\label{NSL-3}
\begin{array}{lcr}
 - \nu \Delta u +\nabla P = - \eta u_t -\nu(\Delta -\Delta_u)u+(\nabla-\nabla_u)P & \mbox{in} & \Omega \times (0,T), \\[4pt]
\div u=\div((\Id-A)u)=-A:Du & \mbox{in} & \Omega \times (0,T), \\[4pt]
u=0 & \mbox{on} & \d \Omega \times (0,T).
\end{array}
\end{equation}

Note that one may use Proposition  \ref{p:bog} so as to handle the potential part 
of $u.$ Therefore using standard results for the stationary Stokes equation (see \cite{Galdi})
enables us to get
\begin{multline}\label{p4}
 \|\nu\nabla^2 u,\nabla P\|_{L_{n^*}(\Omega \times (0,T))} \leq C\big( \|u_t\|_{L_{n^*}(\Omega \times (0,T))} \\
+\|\nu(\Delta-\Delta_u)u,(\nabla-\nabla_u)P,\nu \nabla (A : Du)\|_{L_{n^*}(\Omega \times (0,T))}\big).
\end{multline}
The key to bounding the right-hand side is that, because $n^*>n,$
 we have by embedding and H\"older's inequality
$$
\begin{array}{lll}
 \|\Int_0^t Du(t',y)\,dt'\|_{L_\infty(\Omega \times (0,T))} &\leq& 
 C T^{1\!-\!\frac{1}{n^*}} \|\nabla u\|_{ L_{n^*}(0,T;L_\infty(\Omega))}\\[1ex]
&\leq& CT^{1\!-\!\frac{1}{n^*}} \|\nabla^2 u\|_{ L_{n^*}(\Omega \times (0,T))}.
\end{array}
$$
In particular, this allows to write that 
$$
\begin{array}{lll}
\|(\Delta-\Delta_u)u\|_{L_{n^*}(\Omega \times (0,T))}&\!\!\!\lesssim\!\!\!&
\|D(A{}^T\!A)\|_{L_\infty(0,T;L_{n^*}(\Omega)}\|Du\|_{L_{n^*}(0,T;L_\infty(\Omega)}\\&&\qquad\qquad
+\|\Id-A{}^T\!A\|_{L_\infty(\Omega \times (0,T))}\|D^2u\|_{L_{n^*}(\Omega \times (0,T))},\\[1.5ex]
&\!\!\!\lesssim\!\!\!&
\|D^2u\|_{L_1(0,T;L_{n^*}(\Omega)}\|Du\|_{L_{n^*}(0,T;L_\infty(\Omega))}\\
&&\qquad\qquad+\|Du\|_{L_1(0,T;L_\infty(\Omega))}\|D^2u\|_{L_{n^*}(\Omega\times(0,T))},\\[1.5ex]
&\!\!\!\lesssim\!\!\!&T^{1-1/n^*} \Xi_{(u,P)}^2(T).
\end{array}
$$
Similar estimates hold true for the other terms of the right-hand side of \eqref{p4}.
So finally, putting together all the above inequalities leads to 
\begin{equation}
 \Xi_{(u,P)}(T) \leq C_{u_0}\bigl(1+T^{\delta}\Xi_{(u,P)}(T)\bigr)^3
\end{equation}
for some $\delta=\delta(n)\in(0,1)$ which may be computed explicitly.
 \smallbreak
Then we are able to close the estimate, namely to write that 
\begin{equation}\label{p8}
 \Xi_{(u,P)}(T)\leq 8C_{u_0}
\end{equation}
whenever $T$ has been chosen so that
\begin{equation}\label{eq:time}
8C_{u_0}T^\delta\leq1.
\end{equation}


\subsection{The proof of existence}

In this short subsection, we explain how the proof of existence may be achieved
from the above a priori estimates.

Taking for granted the proof of the existence of a solution in the smooth case is the shortest way. 
Under the assumptions of Theorem \ref{th:exist}, one may for instance smooth out
the initial density $\rho_0$ by convolution by a positive mollifier. 
This provides us with a family of smooth approximate densities $(\rho_0^\ep)_{\ep>0}$
satisfying the same lower and upper bound as $\rho_0.$
Then applying the local and existence and uniqueness statement of e.g. \cite{CK} (bounded domain
case) or \cite{D3} (whole space case), one obtains a family of
solutions $(\rho^\ep,v^\ep,\nabla Q^\ep)$ for System \ref{NSE} with data $(\rho_0^\ep,v_0).$
This family of solutions has the required regularity. 
In addition, the possible blow-up of  $(\rho^\ep,v^\ep,\nabla Q^\ep)$ at time $T$ is controlled
by the norm of $Dv^\ep$ in $L_\infty(0,T;L_2(\Omega))\cap L_1(0,T;L_\infty(\Omega)).$
Note that Proposition \ref{p:change} ensures
that $(\rho^\ep,v^\ep,\nabla Q^\ep)$ corresponds to a solution $(\eta^\ep,u^\ep,\nabla P^\ep)$ of \eqref{NSL}
with the same regularity. 

Now, the computations that have been performed in the previous section, combined with the
aforementioned blow-up criterion ensure that the lifespan of $(\eta^\ep,u^\ep,\nabla P^\ep)$ (or
of $(\rho^\ep,v^\ep,\nabla Q^\ep)$) may be bounded by below as in \eqref{eq:time}, and 
that \eqref{p8} is satisfied. The important point is that all those bounds depend on the density
\emph{only through its infimum and supremum}.
So eventuallly,  $(\rho^\ep,v^\ep,\nabla Q^\ep)$ is uniformly bounded in
$$
L_\infty(\Omega\times(0,T))\times W^{2,1}_{n^*,n^*}(\Omega\times(0,T))
\times L_{n^*}(\Omega\times(0,T)),
$$
and in addition,  $\d_tv^\ep$ is bounded in $L_\infty(0,T;L_2(\Omega))\cap L_2(0,T;W^1_2(\Omega)).$
\medbreak
By resorting to standard compactness argument, it is now 
easy to conclude that this family converges, up to extraction, 
to some $(\rho,v,\nabla Q)$ with the same regularity and satisfying  the same bounds. 
The regularity is so high that the it is clear that it
satisfies \eqref{NSE}. 
 Uniqueness then follows from Theorem \ref{th:uniq}.

\begin{rem}
An alternative approach to the issue of existence can be done by an iterative scheme performed in the same way as in our recent work \cite{DM}, or as in \cite{MZ,Sol} for the homogeneous Navier-Stokes equations in the Lagrangian coordinates.
\end{rem}


\section{Global existence}

This section is dedicated to the proof of \emph{global-in-time} solutions. 
As pointed out in the introduction, in the  case of smooth data, this is a classical 
issue that has been solved by  different authors in the Eulerian framework :
if there is no vacuum initially then global existence may be achieved 
for general (smooth) data in the two-dimensional case and
if the velocity is small in the three-dimensional case (see e.g. \cite{D3,LS}).

As for us, in order to show the global existence,  one may  adopt
 the Eulerian approach, too.  However the very low regularity of the density
 will enforce us to treat the inhomogeneity of the fluid as a perturbation
 (hence to assume \eqref{den-str} or \eqref{den-str-1})
 and to use the Lagrangian framework to prove the uniqueness.

 \subsection{The two-dimensional case}

Here we prove Theorem \ref{th:lar}.
We concentrate  on  the proof of global a priori estimates.
Indeed,  existence   can be established by an elementary approximation with smooth density
exactly as in the previous section : for smooth enough densities,  
the existence  of global solutions with velocity (locally) in  $W^{2,1}_{4,2}(\Omega\times\R_+)$
is ensured by \cite{D} (bounded case) or by \cite{D3} (whole space case). 
In addition, let us  emphasize that, in dimension two, this regularity guarantees that we are allowed to 
change coordinates between the Eulerian and Lagrangian ones.
So the uniqueness follows from  Theorem \ref{th:uniq}.

\medbreak
For getting the global existence, the computations are simpler
in  the Eulerian framework. We aim at  getting a  control 
over $v_t$ in $L_\infty(0,T;L_2(\Omega))$ in terms of the data and of $T$ only. 
Even though it is classical (see e.g. \cite{AKM,D3,LS}) we here
recall how to proceed. 
First, we test the momentum equation of System \eqref{NSE} by $v_t.$ We get: 
 $$
 \int_\Omega\rho|v_t|^2\,dx+\frac\nu2\frac d{dt}\int_\Omega|\nabla v|^2\,dx+\int_\Omega \sqrt\rho v_t\cdot(\sqrt\rho v\cdot\nabla v)\,dx=0.
 $$
 Hence   H\"older and Young inequalities imply that
\begin{equation}\label{H11}
\|\sqrt\rho v_t\|_{L_2(\Omega)}^2+\nu\frac d{dt}\|\nabla v\|_{L_2(\Omega)}^2
\leq \|\sqrt\rho v\|_{L_4(\Omega)}^2
\|\nabla v\|_{L_4(\Omega)}^2.
\end{equation}
On the other hand, using maximal regularity for 
the stationary Stokes equation  
$$\begin{array}{cl}
-\nu\Delta v+\nabla Q=\sqrt\rho\Bigl(\sqrt\rho v_t+\sqrt\rho
v\cdot\nabla v\Bigr)&\qquad\hbox{in }\ \Omega\\
 \div v=0 &\qquad\hbox{in }\ \Omega\\
v=0&\qquad\hbox{on }\ \d\Omega,
\end{array}
$$
gives (omitting the time-dependency)
\begin{equation}\label{H12}
\nu\|\nabla^2 v\|_{L_2(\Omega)}
+\|\nabla Q\|_{L_2(\Omega)}\lesssim  \|\sqrt\rho\|_{L_\infty(\Omega)}
\Bigl(\|\sqrt\rho v_t\|_{L_2(\Omega)}
+\|\sqrt\rho v\|_{L_4(\Omega)}\|\nabla v\|_{L_4(\Omega)}\Bigr).
\end{equation}
Now applying  Ladyzhenskaya inequality 
$\|\nabla v\|_{L_4(\Omega)}^2\lesssim\|\nabla v\|_{L_2(\Omega)}
\|\nabla^2v\|_{L_2(\Omega)},$ 
yields
\begin{eqnarray}\label{H13}
&\nu\|\nabla^2 v\|_{L_2(\Omega)}+\|\nabla Q\|_{L_2(\Omega)}\lesssim
\|\sqrt\rho\|_{L_\infty(\Omega)}\|\sqrt\rho v_t\|_{L_2(\Omega)}\qquad\qquad\qquad\nonumber\\
&\qquad\qquad\qquad\qquad+\Frac{\|\rho\|_{L_\infty(\Omega)}\|\sqrt\rho v\|_{L_4(\Omega)}^2\|\nabla v\|_{L_2(\Omega)}}\nu\cdotp
\end{eqnarray}
 Making use of Ladyzhenskaya inequality  in \eqref{H11}, also leads to 
 \begin{equation}\label{H14}
\|\sqrt\rho v_t\|_{L_2(\Omega)}^2+\nu\frac d{dt}
\|\nabla v\|_{L_2(\Omega)}^2\leq \nu^2\frac{\|\nabla^2 v\|_{L_2(\Omega)}^2}{\|\rho\|_{L_\infty(\Omega)}}
+C\frac{\|\rho\|_{L_\infty(\Omega)}}{\nu^2}\|\sqrt\rho v\|_{L_4(\Omega)}^4\|\nabla v\|_{L_2(\Omega)}^2.
\end{equation}
Finally, adding up  \eqref{H14} and \eqref{H13}, using that 
$\|\rho(t)\|_{L_\infty(\Omega)}=\|\rho_0\|_{L_\infty(\Omega)}$ 
and performing a time integration yields
$$
\displaylines{
\|\nabla v(t)\|_{L_2(\Omega)}^2+\int_0^t\biggl(\frac{\|\sqrt{\rho}v_t\|_{L_2(\Omega)}^2}{2\nu}+\frac{\|\nabla Q\|_{L^2}^2}{\nu\|\rho_0\|_{L_\infty(\Omega)}}
+\nu\frac{\|\nabla^2 v\|_{L_2(\Omega)}^2}{\|\rho_0\|_{L_\infty(\Omega)}}\biggr)\,d\tau
\hfill\cr\hfill\leq \|\nabla v_0\|_{L_2(\Omega)}^2
+\frac{C\|\rho_0\|_{L_\infty(\Omega)}}{\nu^3}\int_0^t\|\sqrt\rho v\|_{L_4(\Omega)}^4\|\nabla v\|_{L_2(\Omega)}^2\,d\tau,}
$$
and Gronwall lemma implies that
\begin{multline}\label{l14}
\nu\|\nabla v(t)\|_{L_2(\Omega)}^2\!+\!\int_0^t\!\bigl(\|\sqrt\rho v_t\|_{L_2(\Omega)}^2
+\|\rho_0\|_{L_\infty(\Omega)}^{-1}\|\nabla Q\|_{L_2(\Omega)}^2
+\nu^2\|\rho_0\|_{L_\infty(\Omega)}^{-1}\|\nabla^2v\|_{L_2(\Omega)}^2\bigr)d\tau\\
\leq \nu \|\nabla v_0\|_{L_2(\Omega)}^2\:e^{C\nu^{-3}\|\rho_0\|_{L_\infty(\Omega)}
\int_0^t\|\sqrt\rho v\|_{L_4(\Omega)}^4\,d\tau}.
\end{multline}

Note that the exponential term is controlled thanks to the basic energy equality \eqref{energy}
(combined with Ladyzhenskaya's inequality).

Since $W^1_2(\R^2)$ is not embedded into $L_\infty(\R^2)$ we  still do not control the change of coordinates
so that we cannot apply Theorem \ref{th:uniq} to get uniqueness.
So we
are required to improve the regularity of the  solution to \eqref{NSE}. 
In fact, it turns out to be possible to obtain $W^{2,1}_{q,p}$ smoothness for any $1<p<\infty$ and $n<q<\infty$ 
via bootstrap method. To avoid technicality, 
we focus on the case $p=2$ and $q=4$ which suffices both to perform the change of coordinates
and to apply the uniqueness result stated in Theorem~\ref{th:uniq}.
We rewrite System \eqref{NSE} as
\begin{equation}\label{l15}
\begin{array}{lcr}
m v_t  -\nu \Delta v + \nabla Q=(m-\rho)v_t -\rho v \cdot \nabla v & \mbox{ in } & \Omega \times (0,T), \\[3pt]
\div v=0  & \mbox{ in } & \Omega \times (0,T), \\
v=0 & \mbox{in}& \d \Omega \times (0,T).
\end{array}
\end{equation}
where $m =\inf_{y\in \Omega} \rho_0(y)$. Note that  the method of characteristics 
ensures that the initial density controls lower and upper pointwise bounds of the density over 
$\Omega \times (0,T)$.

Then using Theorem \ref{th:stokes} we get:
\begin{multline}\label{l16}
\sup_{0\leq t\leq T}\sqrt{m\nu}\|v(t)\|_{B^1_{4,2}(\Omega)}+
\|m v_t, \nu\nabla^2 v,\nabla Q \|_{L_2(0,T;L_4(\Omega))} \\
\leq C\bigl(\|(\rho-m)v_t\|_{L_2(0,T;L_4(\Omega))} +\|\rho v \cdot \nabla v\|_{L_2(0,T;L_4(\Omega))} +\sqrt{m\nu}\|v_0\|_{B^1_{4,2}(\Omega)}\bigr).
\end{multline}
Now, we have
\begin{equation}\label{l17}
\|\rho v \cdot \nabla v\|_{L_2(0,T;L_4(\Omega))} \leq C \|v\|_{L_4(0,T;L_\infty(\Omega))}\|\nabla v\|_{L_4(\Omega\times(0,T))}.
\end{equation}
The right-hand side of \eqref{l16} is bounded by means of \eqref{l14} as
$$W^{2,1}_{2,2} \subset L_4(0,T;L_\infty(\Omega)\cap W^1_4 (\Omega)).$$
The first term of the right-hand side of \eqref{l16} can be absorbed by the left-hand side
  provided $c$ is sufficiently small in \eqref{den-str}. 
This enables us  to justify that the velocity $v$ remains in $W^{2,1}_{4,2}(\Omega\times(0,T))$ for all $T>0.$

Finally, as $4\geq2={\rm dim} \, \Omega$, we are allowed to apply Theorem \ref{th:uniq} in order 
to get  the uniqueness of our constructed solutions. Theorem \ref{th:lar} is proved.


\subsection{Global existence in the  $n$-dimensional case}

In this part we address the global solvability issue in bounded $n$-dimensional 
domains with $n\geq3$. We adopt the Lagrangian framework (however the Eulerian framework 
may be used as well, as regards the existence theory). The result is based on the technique for the homogeneous system
performed in \cite{Mu}.

In contrast with the other sections, working in bounded domains is important: this is due to 
the following result which ensures the exponential decay of the energy norm. 
\begin{lem}
Let $u$ be a sufficiently smooth solution to \eqref{NSL}. Then
\begin{equation}\label{e1}
\frac{1}{2}\frac{d}{dt}\int_\Omega \eta |u|^2 \,dy +\nu \int_\Omega |\nabla_uu|^2 \,dy =0,
\end{equation}
as long as the Lagrangian coordinates are defined.
In addition if  $\Omega$ is bounded then
\begin{equation}\label{e2}
\|u(t)\|_{L_2(\Omega)}^2 \leq e^{-\frac{\nu\lambda_1}{\eta^*} t}\|u_0\|_{L_2(\Omega)}^2.
\end{equation}
where $\lambda_1$ stands for the first eigenvalue of the Laplace operator, and $\eta^*=\|\eta\|_{L_\infty(\Omega)}.$
\end{lem}
\begin{p}
The proof is similar to that of  Lemma \ref{l:ene} : testing $(\ref{NSL})_2$ by $u$ we get \eqref{e1}. 

In order to get \eqref{e2}, it suffices to notice that, owing to incompressibility, we have
$$
\int_\Omega \eta |u|^2 \,dy=\int_\Omega \rho |v|^2 \,dx\quad\hbox{and}\quad
\int_\Omega |\nabla_uu|^2\,dy=\int_\Omega |\nabla v|^2\,dx.
$$
Hence using \eqref{energy} and Poincar\'e's inequality, 
we readily get \eqref{e2}.
\end{p}

\subsubsection*{Proof of Theorem \ref{th:bdd}:} We focus on the proof of global a priori
estimates for smooth solutions to \eqref{NSL}.
Indeed, from those estimates, it is easy to proceed as in Section \ref{s:existence} so as to prove the 
existence of a global solution under the assumptions of Theorem \ref{th:bdd}: 
this is only a matter of smoothing out the initial density so as to construct 
a sequence of smooth solutions (given by e.g.  \cite{D}) with uniform norms. 
\smallbreak
So given a global solution $(\eta,u,\nabla P)$ to \eqref{NSL} with data $\rho_0\in L_\infty(\Omega)$ satisfying
\eqref{den-str-1} and $u_0\in B^{2-2/p}_{q,p}(\Omega)$ with $\div u_0=0$ and $u_0|_{\d\Omega}=0,$ 
 we  introduce the following quantities: 
$$\displaylines{
M_{-1}:= m^{1/p}\nu^{1/p'} \|u_0\|_{B^{2-2/p}_{q,p}(\Omega)}+m\|u_0\|_{L_2(\Omega)},\cr
M_k:=  m^{1/p}\nu^{1/p'} \|u\|_{L_\infty(k,k+1;B^{2-2/p}_{q,p}(\Omega))}+  \|mu_t, \nu \nabla^2 u, \nabla P\|_{L_p(k,k+1;L_q(\Omega))}}
$$
where $m:=\inf \rho_0$ and $k\in\N.$ Recall $1<p<\infty$, $n<q<\infty$.

\medbreak

Let us notice that setting 
$$
u(t,x)=\nu \tilde u(\nu t,x)\quad\hbox{and}\quad
P(t,x)=\nu^2\tilde P(\nu t,x)
$$
reduces our study to the case $\nu=1.$
Hence we shall assume from now on that $\nu=1.$
\medbreak

Define a smooth function $\zeta: \R \to [0,1]$ such that
\begin{equation}\label{b2}
\zeta^k(t)=\left\{ 
\begin{array}{cl}
1 &\hbox{if }\  t \geq 0, \\
0 & \hbox{if }\  t\leq -1,
\end{array}
\right.
\end{equation}
and set $\zeta^k(t):=\zeta(t-k)$ for $k\geq0,$ and $I_k:=[k-1,k+1]$ for $k\geq1.$
\smallbreak

We recast System \eqref{NSL}  with $t_0=k-1$ as follows:
\begin{equation}\label{b3}
\begin{array}{lcr}
m[\zeta^k u]_t - \Delta [\zeta^ku] +\nabla [\zeta^k P] = \zeta^k(m-\eta)u_t && \\[3pt] \qquad\qquad-(\Delta -\Delta_u)[\zeta^ku] +
 (\nabla -\nabla_u) [\zeta^k P] + m(\zeta^k)_t u & \mbox{in} & \Omega \times (0,T),\\[5pt]
\div [\zeta^ku]=\div [\zeta^ku] -\divu [\zeta^k u]  & \mbox{in} & \Omega \times (0,T),\\[5pt]
\zeta^ku=0 & \mbox{on} & \d \Omega \times (0,T),\\[5pt]
\zeta^ku|_{t=k-1}=0 & \mbox{in} & \Omega.
\end{array}
\end{equation}
Let $m^*:=\sup\rho_0.$
We claim that  there exist two positive constants $K$ and $\alpha$ depending only on $m^*,n,\Omega,p,q,$ so that, under Condition \eqref{den-str-1}, we have
\begin{equation}\label{b4}
M_k \leq  KM_{-1} e^{-\alpha k}\quad\hbox{for all }\ k\in\N.
\end{equation}
Let us observe that, by Sobolev embedding (here we use that $q>n$), we have  
\begin{equation}\label{b11}
 \int_0^{k+1} \|Du\|_{L_\infty(\Omega)}\, ds \leq C \sum_{\ell=0}^{k} \|D^2u\|_{L_p(I_\ell;L_q(\Omega))}\leq
C \sum_{\ell=0}^{k} M_k.
\end{equation}
So, given that 
$$
\sum_{k\geq0}  KM_{-1} e^{-\alpha k}=\frac{KM_{-1}}{1-e^{-\alpha}},
$$
if we  assume that $M_{-1}$
is  small enough --a condition which is equivalent to the smallness of $c'$ in \eqref{den-str-1}--
then  \eqref{eq:smallu} is satisfied on $[0,k+1]$ if \eqref{b4} is satisfied up to $k.$
\medbreak
Proving \eqref{b4} will be done by induction on $k.$
The first step, $k=0,$ is clear. This is a direct consequence of Theorem \ref{th:stokes} 
applied to \eqref{NSL} on the time interval $[0,1],$
and of estimates for $A.$ 
\medbreak
Let us now take for granted Inequality \eqref{b4} up to $k-1.$ In order to prove it for $k,$ 
we shall  estimate $(\zeta^ku,\zeta^kP)$ on  the interval $I_k.$
For that,  one may resort once again  to Theorem \ref{th:stokes}. 
First, we bound the right-hand side  of \eqref{b3} in $L_p(I_k;L_q(\Omega))$: we readily have
$$
\begin{array}{lll}
\|\zeta^k(m-\eta) u_t\|_{L_p(I_k;L_q(\Omega))} &\!\!\!\!\leq\!\!\!\!&  (m^*-m)\|u_t\|_{L_p(I_k;L_q(\Omega))},\\[1.5ex]
 \|(\Delta -\Delta_{u})[\zeta^ku]\|_{L_p(I_k;L_q(\Omega))} &\!\!\!\!\leq\!\!\!\!& \|A{}^T\!A-\Id\|_{L_\infty(\Omega\times I_k)}
 \|\nabla^2u\|_{L_p(I_k;L_q(\Omega))}\\&&\qquad
  +\|\nabla(A{}^T\!A)\|_{L_\infty(I_k;L_q(\Omega))}\|\nabla u\|_{L_p(I_k;L_\infty(\Omega))},\\[1.5ex]
  \|(\nabla -\nabla_{u}) [\zeta^k P]\|_{L_p(I_k;L_q(\Omega))}&\!\!\!\!\leq\!\!\!\!&  \|A-\Id\|_{L_\infty(\Omega\times I_k)}
  \|\nabla P\|_{L_p(I_k;L_q(\Omega))},\\[1.5ex]
  \|m(\zeta^k)_t u\|_{ L_p(I_k;L_q(\Omega))}& \!\!\!\!\leq\!\!\!\!& m^*\|u\|_{ L_p(I_k;L_q(\Omega))}. 
\end{array}
$$
Let us notice that, by interpolation and because $\Omega$ is bounded, we have for some $\theta\in(0,1),$ 
$$
\|u\|_{ L_p(I_k;L_q(\Omega))}\leq C \|D^2u\|_{ L_p(I_k;L_q(\Omega))}^\theta\|u\|_{L_p(I_k;L_2(\Omega))}^{1-\theta}.
$$
Therefore, taking advantage of \eqref{e2} and of the definition of $M_{k-1}$ and of $M_k,$ we get for some $\beta>0$ (depending only
on $\Omega,$ $p,$ $q$, $m_*$ and 
for all $\ep\in(0,1),$ 
\begin{equation}\label{eq:decay}
m^*\|u\|_{ L_p(I_k;L_q(\Omega))}\leq \ep(M_{k-1}+M_k) +C_\ep M_{-1}e^{-\beta k}.
\end{equation}
Next, we have to bound the left-hand side of $\eqref{b3}_2$: we have
$$
\begin{array}{lll}
\|\div((\Id \!-\! A) [\zeta^k u])\|_{ L_p(I_k;W^1_q(\Omega))}&\!\!\!\!\leq\!\!\!\!& \|A-\Id\|_{L_\infty(\Omega\times I_k)}
 \|\nabla^2u\|_{L_p(I_k;L_q(\Omega))}\\&&\qquad
  +\|\nabla A\|_{L_\infty(I_k;L_q(\Omega))}\|\nabla u\|_{L_p(I_k;L_\infty(\Omega))} ,\\[1.5ex]
\|\d_t\left((\Id-A)\zeta^k u\right)\|_{ L_p(I_k;L_q(\Omega))}&\!\!\!\!\leq\!\!\!\!& 
\|A-\Id\|_{L_\infty(\Omega \times (0,T))}\|(\zeta^ku)_t\|_{L_p(I_k;L_q(\Omega))}\\&&\qquad\qquad+\|A_t\|_{L_p(I_k;L_\infty(\Omega))} 
\|u\|_{ L_\infty(I_k;L_q(\Omega))}.
\end{array}
$$
Let us look at the quantities depending on the matrix $A$. Recall that
$$A^{-1}=\Id + \int_{0}^tDu(s)\,ds,$$ so taking advantage of \eqref{b11} 
and of the hypothesis that follows, 
one may write that 
$$\|\Id-A\|_{L_\infty(\Omega\times(0,k+1))} \leq 2\|Du\|_{L_1(0,k+1;L_\infty(\Omega))},$$
and a similar inequality for $\Id-A{}^T\!A.$ 
Likewise, we have 
\begin{equation}\label{b13}
\|DA\|_{L_\infty(I_k;L_q(\Omega))} \leq C\sum_{\ell=0}^{k} \|D^2 u\|_{L_q(I_\ell;L_q(\Omega))} \leq
C\sum_{\ell=0}^{k} M_\ell
\end{equation}
and
\begin{equation}\label{b14a}
\|A_{t}\|_{L_p(I_k;L_\infty(\Omega))} \leq C \|Du\|_{L_p(I_k;L_\infty(\Omega)}\leq C(M_{k-1}+M_k).
\end{equation}
So finally, putting together all the previous inequalities and applying Theorem \ref{th:stokes}
to \eqref{b3}, we end up with 
$$\displaylines{
  \|m(\zeta^ku_t), \nabla^2(\zeta^ku), \nabla(\zeta^kP)\|_{L_p(I_k;L_q(\Omega))} \hfill\cr\hfill\leq 
  C(M_{k-1}+M_k) \biggl(\Bigl(\frac{m^*-m}{m}\Bigr)+  \Bigl(\sum_{\ell=0}^{k} M_\ell\Bigr) +\ep\biggr)
 +C_\ep M_{-1}e^{-\beta k}.}
  $$
At this point, it is clear that one has to take $\alpha=\beta.$  Note also that if $M_{-1}$   
and the oscillations of the density are small enough then, taking  $\ep$ small enough too, 
the above inequality implies, up to a change of $C,$ 
$$
M_k\leq C(c M_{k-1} + M_{-1}e^{-\alpha k}),
$$
\medbreak
Now, using the induction hypothesis \eqref{b4} for $M_{k-1},$ we deduce that
$$
M_k\leq KM_{-1}e^{-\alpha k}\biggl(\frac{cCe^{\alpha}}K+\frac CK\biggr).
$$
Therefore, we see that if we take $K=2C$ and assume that $c$ has been chosen so that $cCe^{\alpha}\leq1/2$
then we get \eqref{b4} for $M_k.$
\medbreak
Note that our  proof is not quite rigorous  as 
we did use \eqref{b4} at rank $k$ in the above inequalities. 
To make the argument work, it is just a matter of replacing the interval 
$I_k$ with $[k-1,T].$  By continuity of the norms with respect to time, it is clear that
\eqref{b4} at rang $\ell\leq k-1$ ensures that the desired inequality is satisfied on $[k,T]$ for any $T$
close enough to $k.$ Then resorting to a standard bootstrap argument allows to conclude to the desired
inequality for $M_k.$ 
This completes the proof of  Theorem \ref{th:bdd}.


\section{Appendix}

Throughout this paper, we used repeatedly the following
well-known result 
for the divergence equation (see e.g. \cite{Galdi} and the references therein):
\begin{lem}\label{l:bogovskii}   Let $\Omega$ be a bounded Lipschitz domain of $\R^n.$
There exists a linear  operator $\cB$ which is bounded from $L_q(\Omega)$ to $W^1_q(\Omega)$ for all  $q\in(1,\infty)$ and 
such that for any
$f\in L_q(\Omega)$ the vector-field
$u:=\cB(f)$  satisfies 
\begin{equation}\label{eq:div0}
\div u=f \quad\hbox{in }\ \Omega\qquad\hbox{and}\qquad
u|_{\d\Omega}=0\quad\hbox{on }\ \d\Omega.
\end{equation}
\end{lem}
This result may be proved by means  of an explicit formula -- the \emph{Bogovski\u{\i} formula} -- that provides a solution to the above
divergence equation  in   the case where $\Omega$ is  star-shaped.
In our paper, we had to use a more elaborate version of the above lemma, namely the following statement that has 
been established in \cite{DM-luminy}:
\begin{prop}\label{p:bog}
   Let   $\Omega$ be a $C^2$ bounded domain.
          There exists a linear operator  $B$  acting on couples $(R,\zeta)$ with 
          $R:\Omega\rightarrow\R^n$ and $\zeta:\d\Omega\rightarrow\R$ which is 
          continuous from $L_q(\Omega;\R^n)\times  W^{-1/q}_q(\d \Omega,\R)$ to $L_q(\Omega,\R^n)$ for
          all $q\in(1,+\infty)$ and such that  $u:=B(R,\zeta)$ satisfies the generalized divergence equation:
          \begin{equation}\label{eq:gdiv}
          -\int_\Omega u \cdot \nabla \phi\, dx= -\int_\Omega R \cdot \nabla \phi\, dx
  + \int_{\d \Omega} \zeta\phi\, d\sigma \ 
 \hbox{ for all }\  \phi \in C^\infty(\overline{\Omega}).
  \end{equation}
    If in addition $\div R\in L_q(\Omega)$ and $R\cdot\vec n=0$   then $u:=B(R,0)$ satisfies
    $$
    \div u=\div R\quad\hbox{in }\ \Omega\qquad\hbox{and}\qquad
u|_{\d\Omega}=0\quad\hbox{on }\ \d\Omega,
$$
  and the following inequality holds true:
   \begin{equation}\label{eq:div1}
     \|u\|_{W^{1}_q(\Omega)} \leq C\|\div R\|_{L_{q}(\Omega)}.
     \end{equation}
     Furthermore, if we also have $\div R\in W_q^1(\Omega)$ then $u$ is in $W^2_q(\Omega)$ and we have
   \begin{equation}\label{eq:div2}
     \|u\|_{W^{2}_q(\Omega)} \leq C\|\div R\|_{W^1_{q}(\Omega)}.
     \end{equation}
\end{prop}

We claim that this  statement implies Lemma \ref{l:bog}. Indeed, we set 
$u:=B(R,0).$ Then it is is clear that  \eqref{eq:bog1} holds true. 
Then differentiating $u$ with respect to time yields 
$$
u_t=B(R_t,0).
$$
Hence applying the first part of the above statement yields \eqref{eq:bog2}.
\bigbreak
In Section \ref{s:existence}, owing to our use of Lagrangian coordinates, 
it was natural to extend Lemma \ref{l:bog} and Proposition \ref{p:bog}
to the \emph{twisted} divergence equation, namely
$$
\divA u=\div R \ \hbox{ in }\ \Omega\qquad\hbox{and }\ 
u=0\ \hbox{ on }\ \d\Omega
$$
with  $\divA u:=\div(A u).$
In particular, we used the following statement:
\begin{lem}\label{l:divA} 
 Let   $\Omega$ be a $C^2$ bounded domain.
  Let $A\in L_\infty (\Omega; \R^{n\times n})$ be such that
$\det A\equiv1.$ 
There exists a positive constant $c$ such that if 
 \begin{equation}\label{eq:smallA}
  \|A-\Id\|_{L_\infty(\Omega)}\leq c
 \end{equation}
 then there exists a map $B_A$ acting on couples $(R,\zeta)$ with 
          $R:\Omega\rightarrow\R^n$ and $\zeta:\d\Omega\rightarrow\R$ which is 
          continuous from $L_q(\Omega;\R^n)\times  W^{-1/q}_q(\d \Omega,\R)$ to $L_q(\Omega,\R^n)$ for
          all $q\in(1,+\infty)$ and such that  $u:=B_A(R,\zeta)$ satisfies the 
          generalized twisted divergence equation:
 \begin{equation}\label{eq:gdivA}
          -\int_\Omega Au \cdot \nabla \phi\, dx= -\int_\Omega R \cdot \nabla \phi\, dx
  + \int_{\d \Omega} \zeta\phi\, d\sigma \ 
 \hbox{ for all }\  \phi \in C^\infty(\overline{\Omega}).
  \end{equation}
    If in addition $\div R\in L_q(\Omega)$ and $R\cdot\vec n=0$   then $u:=B_A(R,0)$ satisfies
\begin{equation}\label{eq:divA}
\divA u=A:Du=\div R \ \hbox{ in }\ \Omega\qquad\hbox{and }\ 
u=0\ \hbox{ on }\ \d\Omega,
\end{equation}
and \eqref{eq:div1}
with a constant independent of $A.$ 
\medbreak
Finally, in the smooth case, if the data $R$ and $A$ depend on a parameter
$t$ in some interval of $\R$ with $R_t$ in $L_q(\Omega)$ and \eqref{eq:smallA} satisfied
for almost all $t,$ then $u$ fulfills:
\begin{equation}\label{eq:ut}
\|u_t\|_{L_q(\Omega)}\leq C\bigl(\|A_tu\|_{L_q(\Omega)}+\|R_t\|_{L_q(\Omega)}\bigr).
\end{equation} 
\end{lem}
\begin{p}
The proof follows from Proposition  \ref{p:bog}: we consider the linear operator $T$ defined by 
\begin{equation} T(\bar \xi)=\xi, \mbox{  where }
\xi := B\bigl((\Id- A) \bar \xi +  R,\zeta\bigr).
\end{equation}
 Note that this definition and the fact that $\det A\equiv1$   imply   that
 any fixed point of $T$ satisfies  \eqref{eq:gdivA} (or \eqref{eq:divA} in the smooth case
 with $\zeta=0$).  
 Next, under Condition \eqref{eq:smallA} with $c$ small enough, 
one may apply the Banach fixed point theorem to $T$ so as to get a solution to our problem.
The reader may refer to \cite{DM-luminy} for more details.
\smallbreak
Concerning the proof of inequality \eqref{eq:ut}, it suffices to differentiate once
the equality
$$
u =B\bigl((\Id- A)u +  R,0\bigr).
$$
We get
$$
u_t=B\bigl(-A_tu+(\Id-A)u_t+R_t,0\bigr).
$$
So it is a mere consequence of the first part of Proposition \ref{p:bog}.
\end{p}

{\bf Acknowledgement}. The second author (PBM) thanks  Vladimir \v{S}verak for a fruitful discussion. The second author has been supported by 
the MN grant IdP2011 000661.

\end{document}